\documentclass{amsart}
\usepackage{amsfonts}

\newtheorem{theorem}{Theorem}
\theoremstyle{plain}

\newtheorem{corollary}{Corollary}

\newtheorem{lemma}{Lemma}

\newtheorem{proposition}{Proposition}

\numberwithin{equation}{section}

\input{tcilatex}

\begin{document}
\title[Stochastic Calculus and Quantum Filtering]{Quantum Stochastic
Calculus and Quantum Nonlinear Filtering}
\author{V. P. Belavkin}
\address{On leave of absence from M.I.E.M., B. Vusovski Street 3/12 Moscow
109028, USSR.}
\address{Centro Matematico V. Volterra Dipartimento di Matematica Universit%
\`{a} di Roma II}
\email{vpb@maths.nott.ac.uk}
\urladdr{http://www.maths.nott.ac.uk/personal/vpb/}
\thanks{Published in: \textit{Journal of Multivariate Analysis}, \textbf{42}
(2) 171--201 (1992).}
\date{20 September, 1989}
\subjclass{}
\keywords{Quantum stochastic calculus, Quantum Langevin equations, Quantum
nondemolition processes, Quantum conditional expectations, Quantum nonlinear
filtering}
\dedicatory{}
\thanks{This paper is in final form and no version of it will be submitted
for publication elsewhere.}

\begin{abstract}
A $\star $--algebraic indefinite structure of quantum stochastic (QS)
calculus is introduced and a continuity property of generalized nonadapted
QS integrals is proved under the natural integrability conditions in an
infinitely dimensional nuclear space. The class of nondemolition output QS
processes in quantum open systems is characterized in terms of the QS
calculus, and the problem of QS nonlinear filtering with respect to
nondemolition continuous measurments is investigated. The stochastic
calculus of a posteriori conditional expectations in quantum observed
systems is developed and a general quantum filtering stochastic equation for
a QS process is derived. An application to the description of the
spontaneous collapse of the quantum spin under continuous observation is
given.
\end{abstract}

\maketitle

\section*{Introduction}

The problem of description of continuous observation and filtering in
quantum dynamical systems can be effectively solved in the framework of
quantum stochastic (QS) calculus of nondemolition input-output processes
first developed for quantum unitary Markovian evolutions in \cite{bib:12}.

In contrast to classical probability theory, the conditional expectations
defining \`a posteriori states of a quantum system with respect to a
subalgebra may not exist in general and the existence depends on the algebra
of observables and on the \`a priori state.

During the preparation of the measurement of a Quantum System the change
necessary to produce the \`a priori compatible state as a mixture of the \`a
posteriori states is referred to in quantum physics as the \textit{demolition%
} of the system. The latter involves a change in the initial state by the
reduction of the algebra of the system during such a preparation. The
nondemolition principle provides a sufficient condition for the algebra of
the quantum system to be prepared for the measurement in an initial state.

The mathematical formulation of the nondemolition principle for the
observability of a class of quantum processes was given in \cite{bib:4} and
investigated in subsequent papers \cite{bib:5,bib:3}. This fundamental
principle of quantum measurement theory means that if a QS process $X_{t}$
is indirectly observable by the measurement of another process $Y_{t}$ then $%
X$ and $Y$ must satisfy the one sided commutativity condition $%
[X_{t},Y_{s}]\equiv X_{t}Y_{s}-Y_{s}X_{t}=0$ for all $t\geq s$ but not for $%
t<s$.

In the physical language this means that the measurements of $Y$ in real
time do not demolish the quantum system $X$ (which has been prepared for the
observation) at the present time or in the future. The condition given
above, however, shows that, though the past of $X$ (priori to $t$) can never
be observed, it is demolished by the observation of the $Y$ process.
Mathematically it can be expressed as the decomposability of the algebra $%
\mathcal{A}_{t}$ generated by $\{X(s):s\geq t\}$, describing the present and
future of the system with respect to the spectral resolution of any
Hermitian operator of the algebra $\mathcal{B}^{t}$ generated by $%
\{Y(s):s\leq t\}$.

In this paper we show using method of quantum filtering that the
nondemolition condition given above is necessary and sufficient for the
evaluation of a posteriori mean values of $X\in\mathcal{A}_t$ given are
arbitrary initial state. In other words we prove that a quantum system is
statistically predictable by a measurement procedure, iff the observable
process satisfies the nondemolition condition.

In the Sections \ref{sec:belat1} and \ref{sec:belat2} of the paper we
develop the general QS calculus of nondemolition input--output quantum
processes in Fock space, tensored by an initial Hilbert space. We introduce
the QS calculus of such processes using the $\star $--algebraic Minkowski
metric structure of the basic quantum processes and the simple and
convenient notation developed in \cite{bib:ref16}. The Fock representation
of this structure is closely connected with the Lindsay--Maassen kernel
calculus of \cite{bib:13} but is given in terms of the matrix elements of
operators for general quantum noise in Fock space instead of their kernels.
We define the QS integrals in the framework of the new noncommutative
stochastic analysis in the Fock scale which is described in the first
section.

In the Sections \ref{sec:belat3} and \ref{sec:belat4} we give complete
proofs of the results, first formulated in \cite{bib:ref13}, for the general
(non Markovian) quantum filtering from the viewpoint of QS calculus. The
advantage of the $\star $ -- matrix notation enables us to prove the main
filtering theorem for general output process as by using the indefinite
metric for the corresponding $\star $ -- algebra of generators of these
nondemolition processes.

The Markovian nonlinear filtering problem in the framework of quantum
operational (non-stochastic) approach was first investigated in \cite{bib:8}%
, and the possibility of deriving the stochastic equations of quantum
filtering within this framework was shown in \cite{bib:bel11}. The Markovian
filtering for the quantum Gaussian case and the corresponding quantum Kalman
linear filter, first obtained for the one-dimensional case in \cite%
{bib:4,bib:5}, is considered using the QS calculus approach in \cite%
{bib:ref29}.

The present paper is devoted essentially to the study of the nonlinear
problem, extending the innovation martingale methods of the classical
filtering theory \cite{bib:belat1,bib:belat2} to the noncommutative set up
of our problem. An application of the quantum filtering theory to the
solution of the problem of the continuous observation of quantum spin states
is given in Section \ref{sec:belat5}.

\textbf{Acknowledgements.} I wish to thank Prof. L. Accardi, A. Barchielli,
G. Kallianpur, G. Lupieri and M. Piccioni for stimulating discussions and
useful suggestions during the preparing of the paper. The first part of this
paper was written in the Physics Department of the University of Milan and
the second part in Centro Matematico V. Volterra of the University of Roma
II, for the hospitality of which I am very grateful. 

\section{QS calculus of input Bose processes in Fock space}

\label{sec:belat1}

Let us denote by $\mathcal{F}=\Gamma (\mathcal{E})$ the state space of the
one--dimensional Bose--noise, that is the Fock space over the Hilbert space $%
\mathcal{E}=L^{2}(\mathbb{R}^{+})$ of square--integrable complex functions $%
t\mapsto \varphi (t)$ on the real half--line $\mathbb{R}^{+}$. One should
consider $\mathcal{F}$ as the Hilbert space $\Gamma (\mathcal{E}%
)=L^{2}(\Omega (\mathbb{R}^{+}))$ of the square--integrable functions $\tau
\mapsto \varphi (\tau )$ of $\tau =\left( t_{1},\dots ,t_{n}\right) $ with $%
t_{i}\in \mathbb{R}^{+}$,$\;t_{1}<\dots <t_{n}$,$\;n=0,1,2,\dots $ and
scalar product $<\varphi |\chi >=\int \varphi (\tau )^{\ast }\chi (\tau )%
\mathrm{d}\tau $, 
\begin{equation*}
\int \varphi (\tau )^{\ast }\chi (\tau )\mathrm{d}\tau =\sum_{n=0}^{\infty
}\int_{t_{1}\leq \dots \leq t_{n}}\int {\overline{\varphi }}\left(
t_{1},\dots ,t_{n}\right) \chi \left( t_{1},\dots ,t_{n}\right) \mathrm{d}%
t_{1}\dots \mathrm{d}t_{n},
\end{equation*}%
where the integral is taken over the set $\Omega (\mathbb{R}^{+})$ of all
finite chains $\tau $ on $\mathbb{R}^{+}$ with respect to the natural
Lebesgue measure $\mathrm{d}\tau =\mathrm{d}t_{1}\dots \mathrm{d}t_{n}$ for
every $n=|\tau |=0,1,\dots $. Following \cite{bib:ref16} we shall identify
the chains $\tau =\left( t_{1},\dots ,t_{n}\right) $ with the finite subsets 
$\left\{ t_{1},\dots ,t_{n}\right\} \subset \mathbb{R}^{+}$, so that the
empty chain $(n=0)$ is identified with the empty subset $\tau =\emptyset $
having $\mathrm{d}\tau =1$ and $\tau =t\;(n=1)$ is identified with the
one--point subset $\{t\}$ having $\mathrm{d}\tau =\mathrm{d}t$. We shall
also denote the normalized vacuum function as the Kronecker $\delta $%
-function: $\delta _{\emptyset }(\tau )=1$, if $\tau =\emptyset ;\,\,\delta
_{\emptyset }(\tau )=0$, if $\tau \not=\emptyset $, and consider the Fock
spaces $\mathcal{F}_{s}=\Gamma (\mathcal{E}_{s})$, $\,\mathcal{F}%
_{t}^{s}=\Gamma (\mathcal{E}_{t}^{s})$ over orthogonal subspaces $\mathcal{E}%
_{s}=\{\varphi (t)=0:t\leq s\}$,\thinspace $\mathcal{E}_{t}^{s}=\{\varphi
(r)=0:r\notin \lbrack t,s]\}$ as the function Hilbert spaces $L^{2}(\Omega
_{s})$, $L^{2}(\Omega _{t}^{s})$ on the subsets $\Omega _{s}=\{\tau \subset
]s,\infty \lbrack \}$, $\,\Omega _{t}^{s}=\{\tau \subset ]t,s]\}$ of the
chains $\tau >s,\,\,s\geq \tau >t$ correspondingly.

Note that for any $t>s$ a chain $\tau \in \Omega $ can be represented as the
triple $\tau =\left( \tau ^{t},\tau _{t}^{s},\tau _{s}\right) $ of the
subchains $\tau _{s}=\left\{ t_{i}\in \tau :t_{i}>s\right\} $,$\,\ \,\tau
_{t}^{s}=\left\{ t_{i}\in \tau :s\geq t_{i}>t\right\} $, $\,\tau
^{t}=\{t_{i}\in \tau :t_{i}\leq t\}$ so that the direct product
representation $\Omega =\Omega ^{t}\times \Omega _{t}^{s}\times \Omega _{s}$
holds and, hence, the tensor representation $\mathcal{F}=\mathcal{F}%
^{t}\otimes \mathcal{F}_{t}^{s}\otimes \mathcal{F}_{s}$ with $\mathcal{F}%
^{t}=L^{2}\left( \Omega ^{t}\right) $, $\Omega ^{t}=\{\tau \in {\Omega }%
:\tau \leq t\}$.

The basic processes for QS calculus in Fock space $\mathcal{F}$ are the
annihilation $A_{-}$, creation $A^{+}$ and quantum number $N$ processes,
represented for all $t>0$ by the unbounded operators 
\begin{eqnarray*}
(A_{-}(t)\varphi )(\tau ) &=&\int_{0}^{t}\varphi (\tau \sqcup s)\mathrm{d}s,
\\
(A^{+}(t)\varphi )(\tau ) &=&\sum_{s\in \tau }\chi ^{t}(s)\varphi (\tau
\backslash s),
\end{eqnarray*}%
with the common dense domain $\mathcal{D}^{t}=\left\{ \varphi \in \mathcal{F}%
:\int |\tau ^{t}||\varphi (\tau )|^{2}\mathrm{d}t<\infty \right\} ,$ and 
\begin{equation}
(N(t)\varphi )(\tau )=|\tau ^{t}|\varphi (\tau ),  \label{eq:1.1}
\end{equation}%
where $\chi ^{t}(s)=1$, if $s\leq t,\,\,\chi ^{t}(s)=0$, if $s>t,\,\,|\tau
^{t}|=\sum_{s\in \tau }\chi ^{t}(s)$, the chain $\tau \sqcup s$ is defined
almost everywhere as $\left( t_{1},\dots t_{i},s,t_{i+1},\dots ,t_{n}\right) 
$, if $t_{i}<s<t_{i+1}$, and $\tau \backslash s=\left( t_{1},\dots
,t_{i-1},t_{i+1},\dots t_{n}\right) $, if $s=t$, $\tau \backslash s=\tau $,
if $s\not=t_{i}$ for all $i$. Note that the processes $A_{-},\;A^{+}$ and $N$
are non-commuting, but commuting with increments: 
\begin{eqnarray}
\left[ A_{-}(t),A^{+}(t^{\prime })\right] &=&t\wedge t^{\prime }I,\quad 
\mathrm{where}\,\,t\wedge t^{\prime }=\mathrm{min}(t,t^{\prime }),
\label{eq:1.2} \\
\lbrack A_{-}(t),N(t^{\prime })] &=&A_{-}(t\wedge t^{\prime
}),[N(t),A^{+}(t)]=A^{+}(t\wedge t^{\prime }),  \notag
\end{eqnarray}%
the processes $A_{-}$ and $A^{+}$ are mutually adjoint: $A_{-}^{\ast
}(t)=A_{-}(t)^{\ast }=A^{+}(t)$, and $N$ is selfadjoint: $N^{\ast }=N$.

Let us introduce the notations \cite{bib:ref16} 
\begin{equation}
A_{-}^{+}(t)=tI,\,A_{-}^{\circ }(t)=A_{-}(t),\,A_{\circ
}^{+}=A^{+}(t),A_{\circ }^{\circ }(t)=N(t),  \label{eq:1.3}
\end{equation}%
where $I$ is the identity operator in $\mathcal{F}$, thus defining a $%
3\times 3$ matrix--valued QS process $\mathbf{A}=\left( A_{\nu }^{\mu
}\right) $, indexed by $\mu ,\nu \in \{-,o,+\}$ with $A_{\nu }^{\mu }=0$, if 
$\mu =+$ or $\nu =-$. We shall consider the process $\mathbf{A}$ defined as
a linear operator-valued function $A(\mathbf{c},t)=\mathrm{tr}\{\mathbf{c}\;%
\mathbf{A}(t)\}$ in terms of a $3\times 3$ -- matrix $\mathbf{c}=(c_{\nu
}^{\mu })$, 
\begin{equation}
A(\mathbf{c},t)=Ic_{+}^{-}t+A_{-}(c_{\circ }^{-},t)+A^{+}(c_{+}^{\circ
},t)+N(c_{\circ }^{\circ },t),  \label{eq:1.4}
\end{equation}%
where $A_{-}(c_{\circ }^{-})=c_{\circ }^{-}A_{-},A^{+}(c_{+}^{\circ
})=c_{+}^{\circ }A^{+},\,\,N(c_{\circ }^{\circ })=c_{\circ }^{\circ }N$,
writing the matrix trace as $\mathrm{tr}\{\mathbf{cA}\}=c_{\nu }^{\mu
}A_{\mu }^{\nu }$ by the tensor notation of the sum $\sum c_{\nu }^{\mu
}A_{\mu }^{\nu }$. The matrices $\mathbf{c}$ with $c_{\nu }^{\mu }=0$ for $%
\mu =+$ or $\nu =-$ form a complex Lie $\star $-algebra with respect to the
matrix commutator and the involution 
\begin{equation}
\mathbf{c}=\left( 
\begin{array}{ccc}
0 & c_{\circ }^{-} & c_{+}^{-} \\ 
0 & c_{\circ }^{\circ } & c_{+}^{\circ } \\ 
0 & 0 & 0%
\end{array}%
\right) \mapsto \mathbf{c}^{\star }=\left( 
\begin{array}{ccc}
0 & c_{+}^{\circ \ast } & c_{+}^{-\ast } \\ 
0 & c_{\circ }^{\circ \ast } & c_{\circ }^{-\ast } \\ 
0 & 0 & 0%
\end{array}%
\right) =\mathbf{gc}^{\dag }\mathbf{g},\ \mathbf{g}=\left( 
\begin{array}{ccc}
0 & 0 & 1 \\ 
0 & 1 & 0 \\ 
1 & 0 & 0%
\end{array}%
\right) \ ,  \label{eq:1.5}
\end{equation}%
where $c_{\nu }^{\mu \ast }=c_{\nu }^{-\mu }=c_{\mu }^{\ast \nu }$ and $%
\mathbf{g}^{\dag }=\mathbf{g}=\mathbf{g}^{-1}$ is the indefinite metric
matrix, defining a pseudo-scalar product in $\mathbb{C}^{3}$: 
\begin{equation*}
(\mathbf{x}|\mathbf{z})=\overline{x}^{+}z^{-}+\overline{x}^{\circ }z^{\circ
}+\overline{x}^{-}z^{+}=\mathbf{x}^{\star }\mathbf{z}\ ,
\end{equation*}%
$\mathbf{x}^{\star }=(\overline{x}_{+},\overline{x}_{\circ },\overline{x}%
_{-})=\mathbf{x}^{\dag }\mathbf{g}$ is the row, conjugate to the column $%
\mathbf{x}=(x^{\mu })\in \mathbb{C}^{3}$.

Now we can consider a multi-dimensional Bose noise, when $\mathcal{E}$ is a
Hilbert space $L^{2}(${$\mathbb{R}$}$^{+}\rightarrow \mathbb{C}^{m})$ of
vector-functions $\varphi (t)=(\varphi ^{j})(t)\equiv \varphi ^{\circ
}(t),j=1,\dots ,m$ with%
\begin{equation*}
<\varphi |\varphi >=\int \sum_{j=1}^{m}\overline{\varphi }^{j}\varphi ^{j}%
\mathrm{d}t.
\end{equation*}%
It is enough to regard $c_{\circ }^{-}$ as a $m$-row with components $%
c_{j}^{-}\in \mathbb{C},c_{+}^{\circ }$ as a $m$-column with components $%
c_{+}^{j}\in ${$\mathbb{C}$}, and $c_{\circ }^{\circ }$ as a $m\times m$%
-matrix with elements $c_{k}^{i}\in \mathbb{C}$. The following theorems are
valid also for the general situation $\mathcal{F}=\Gamma (L^{2}(\mathbb{R}%
^{+}\rightarrow \mathcal{K}))$, if the indices $\mu ,\nu $ take values in
the set ${-,J,+}$, where the one-point index value $\mu ,\nu =0$ is split
into $m=|J|$ points $j\in J$ of an index set $J$ for a basis in a Hilbert
space $\mathcal{K}$ with the infinite cardinality $|J|=\dim \mathcal{K}$.
\medskip

\begin{proposition}
The basic QS process $A(\mathbf{c})$, defined by (\ref{eq:1.1}), (\ref%
{eq:1.2}), gives for each $t$ an operator representation of the complex Lie $%
\star $-algebra of matrices (\ref{eq:1.5}): $A(\mathbf{c},t)^{\ast }=A(%
\mathbf{c}^{\star },t)$, 
\begin{equation}
\lbrack A(\mathbf{c}^{\star },t),A(\mathbf{c},t^{\prime })]=A([\mathbf{c}%
^{\star },\mathbf{c}],t\wedge t^{\prime })\ .  \label{eq:1.6}
\end{equation}%
The multiplication table \cite{bib:12} for Ito differentials \textrm{d}$%
A_{-},\mathrm{d}A^{+},\mathrm{d}N$, and $I\mathrm{d}t$ can be written in
terms of $A(\mathbf{c},\mathrm{d}t)=c_{\nu }^{\mu }\mathrm{d}A_{\mu }^{\nu
}(t)$ as 
\begin{equation}
A(\mathbf{c}^{\star },\mathrm{d}t)A(\mathbf{c},\mathrm{d}t^{\prime })=A(%
\mathbf{c}^{\star }\mathbf{c},\mathrm{d}t\cap \mathrm{d}t^{\prime }),
\label{eq:1.7}
\end{equation}%
where $\mathrm{d}t\bigcap \mathrm{d}t^{\prime }=\emptyset $ for $%
t\not=t^{\prime },A(\cdot ,\emptyset )=0$, and $\mathrm{d}t\bigcap \mathrm{d}%
t^{\prime }=\mathrm{d}t$ for $t=t^{\prime }$. \label{prop:1}
\end{proposition}

\noindent \textsc{Proof.}\ Taking into account, that $A_{-}^{\ast }=A^{+}$
and $N^{\ast }=N$, one obtains 
\begin{equation}
A(\mathbf{c},t)^{\ast }=Ic_{+}^{-\ast }t+A_{-}(c_{+}^{\circ \ast
},t)+A^{+}(c_{\circ }^{-\ast },t)+N(c_{\circ }^{\circ \ast },t)
\label{eq:1.8}
\end{equation}%
The comparing of (\ref{eq:1.6}) with (\ref{eq:1.2}) gives the $\star $%
-property $A(\mathbf{c})^{\ast }=A(\mathbf{c}^{\star })$ of the map $\mathbf{%
c}\mapsto A(\mathbf{c})$.

The Lie representation property follows directly from the canonical
commutation relations 
\begin{eqnarray*}
&[A_{-}(b_{\circ }^{-},t),A^{+}(d_{+}^{\circ },t)]=tb_{\circ
}^{-}d_{+}^{\circ },\;&[N(b_{\circ }^{\circ }),N(d_{\circ }^{\circ
})]=N([b_{\circ }^{\circ },d_{\circ }^{\circ }])\ , \\
&[N(b_{\circ }^{\circ }),A^{+}(d_{+}^{\circ })]=A^{+}(b_{\circ }^{\circ
}d_{+}^{\circ }),\;&[A_{-}(b_{\circ }^{-}),N(d_{\circ }^{\circ
})]=A_{-}(b_{\circ }^{-}d_{\circ }^{\circ })\ ,
\end{eqnarray*}%
which give $[A(\mathbf{b}),A(\mathbf{d})]=A([\mathbf{b},\mathbf{d}])$, where
we take into account that 
\begin{equation*}
(\mathbf{b}\;\mathbf{d})_{\circ }^{-}=b_{\circ }^{-}d_{+}^{\circ },(\mathbf{b%
}\;\mathbf{d})_{\circ }^{-}=b_{\circ }^{-}d_{\circ }^{\circ },(\mathbf{b}\;%
\mathbf{d})_{+}^{\circ }=b_{\circ }^{\circ }d_{+}^{\circ },\;(\mathbf{b}\;%
\mathbf{d})_{\circ }^{\circ }=b_{\circ }^{\circ }d_{\circ }^{\circ }
\end{equation*}%
for matrices $\mathbf{b},\mathbf{d}$ of the form (\ref{eq:1.5}).

Applying it to $\mathbf{b}=\mathbf{c}^{\star },\mathbf{d}=\mathbf{c}$ and
taking into account the commutativity of $A(\mathbf{c}^{\star },t)$ with
increment $A(\mathbf{c},t^{\prime })-A(\mathbf{c},t)$, one obtains (\ref%
{eq:1.6}). In the same way one obtains (\ref{eq:1.7}) from the Hudson --
Parthasarathy multiplication table 
\begin{eqnarray}
&\mathrm{d}A_{-}(b_{\circ }^{-})\mathrm{d}A^{+}(d_{+}^{\circ })=I\mathrm{d}%
t(b_{\circ }^{-}d_{+}^{\circ }),\;&\mathrm{d}A_{-}(b_{\circ }^{-})\;\mathrm{d%
}N(d_{\circ }^{\circ })=\mathrm{d}A_{-}(b_{\circ }^{-}d_{\circ }^{\circ }) 
\notag \\
&\mathrm{d}N(b_{\circ }^{\circ })\mathrm{d}A^{+}(d_{+}^{\circ })=\mathrm{d}%
A^{+}(b_{\circ }^{\circ }d_{+}^{\circ }),\;&\mathrm{d}N(b_{\circ }^{\circ })%
\mathrm{d}N(d_{\circ }^{\circ })=\mathrm{d}N(b_{\circ }^{\circ }d_{\circ
}^{\circ }),  \label{eq:1.9}
\end{eqnarray}%
for $\mathrm{d}A_{\nu }^{\mu }(t)=A_{\nu }^{\mu }(t+\mathrm{d}t)-A_{\nu
}^{\mu }(t),\mathbf{b}=\mathbf{c}^{\star },\mathbf{d}=\mathbf{c}$. Due to
complex linearity of the map $\mathbf{c}\mapsto A(\mathbf{c})$ the formulas (%
\ref{eq:1.6}), (\ref{eq:1.7}) can be always extended to arbitrary $\mathbf{b}%
,\mathbf{d}$ by polarization formula 
\begin{eqnarray*}
A(\mathbf{b}\;\mathbf{d}) &=&\sum_{n=0}^{3}A\left( (\mathbf{b}^{\star }+%
\mathrm{i}^{n}\mathbf{d})^{\star }(\mathbf{b}^{\star }+\mathrm{i}^{n}\mathbf{%
d})\right) /4\mathrm{i}^{n}\quad ,\quad \mathrm{i}=\sqrt{-1}, \\
A(\mathbf{b})A(\mathbf{d}) &=&\sum_{n=0}^{3}A(\mathbf{b}^{\star }+\mathrm{i}%
^{n}\mathbf{d})^{\star }A(\mathbf{b}^{\star }+\mathrm{i}^{n}\mathbf{d})/4%
\mathrm{i}^{n}.
\end{eqnarray*}%
Hence, (\ref{eq:1.6}) is equivalent to (\ref{eq:1.8}) and (\ref{eq:1.7}) to (%
\ref{eq:1.9}).\hfill $\vrule height.9exwidth.8exdepth-.1ex$

Let us now define a QS integral with respect to the basic process $A$ for a
matrix quantum process $\mathbf{C}(t)=\left( C_{\nu }^{\mu }\right) (t)$, $%
\mu ,\nu \in \{-,J,+\}$ in $\mathcal{F}$. Assuming that the operator--valued
functions $t\mapsto C_{\nu }^{\mu }(t)$ are weakly measurable and adapted: $%
C(t)=C^{t}\otimes I_{t}$, where $C^{t}$ are the operators in $\mathcal{F}%
^{t} $ for all $\mu \in \{-,J\}$ and $\nu \in \{J,+\}$, one can define in
the case of finite $J$ the QS--integral%
\begin{equation*}
\int_{0}^{t}A(\mathbf{C},\mathrm{d}s):=\int_{0}^{t}\sum_{\mu ,\nu }C_{\nu
}^{\mu }\mathrm{d}A_{\mu }^{\nu }\equiv \int_{0}^{t}C_{\nu }^{\mu }\mathrm{d}%
A_{\mu }^{\nu }
\end{equation*}
as the sum of the Lebesgue operator--valued integral $\int C_{+}^{-}(s)%
\mathrm{d}s$ and the It\^{o} integrals $\int C_{j}^{-}\mathrm{d}A_{-}^{j}$,$%
\;\int C_{+}^{j}\mathrm{d}A_{j}^{+}$,$\;\int C_{k}^{i}\mathrm{d}N_{i}^{k}$
in the sense \cite{bib:belat3,bib:lis26}.

In the general case $\mathcal{E}=L^{2}(\mathbb{R}^{+}\rightarrow \mathcal{K}%
) $ we shall regard the QS--integral $\int C_{\nu }^{\mu }\mathrm{d}A_{\mu
}^{\nu }$ as a continuous operator $\mathcal{F}^{+}\rightarrow \mathcal{F}%
_{-}$ on the projective limit $\mathcal{F}^{+}=\bigcap_{\eta >1}\mathcal{F}%
(\eta )$ into $\mathcal{F}_{-}=\bigcap_{\eta <1}\mathcal{F}(\eta )$ of
Hilbert spaces $\mathcal{G}(\zeta )\subset \mathcal{F}\subset \mathcal{G}%
(\xi )$, $\zeta >1>\xi $, with respect to the scalar products 
\begin{equation*}
\Vert \varphi \Vert ^{2}(\eta )=\int_{\eta }^{|\tau |}\eta ^{|\tau |}%
\begin{array}{c}
\langle \varphi (\tau )|\varphi (\tau )\rangle _{+} \\ 
\langle \varphi (\tau )|\varphi (\tau )\rangle _{-}%
\end{array}%
\,\mathrm{d}\tau \ ,\quad 
\begin{array}{c}
\varphi (\tau )\in \mathcal{E}(\tau )\ ,\eta >1 \\ 
\varphi (\tau )\in \mathcal{E}^{\prime }(\tau )\ ,\eta <1\ .%
\end{array}%
\end{equation*}%
Here $\langle \varphi |\varphi \rangle _{+}(\tau )\geq \Vert \varphi \Vert
^{2}\geq \langle \varphi |\varphi \rangle _{-}(\tau )$ are the square-norms
in the Hilbert tensor products $\mathcal{E}(\tau )=\otimes _{t\in T}\mathcal{%
E}(t)$, $\mathcal{K}^{\otimes |\tau |}$, $\mathcal{E}^{\prime }(\tau
)=\otimes _{t\in T}\mathcal{E}^{\prime }(t)$ of Hilbert spaces $\mathcal{E}%
(t)\subseteq \mathcal{K}\subseteq \mathcal{E}^{\prime }(t)$, forming a
Gelfand triple for each $t\in \mathbb{R}^{+}$ with respect to the scalar
product $\Vert \varphi \Vert ^{2}=\langle {\varphi }|{\varphi }\rangle $ in
a Hilbert space $\mathcal{K}$ (or simply $\mathcal{E}(t)=\mathcal{K}=%
\mathcal{E}^{\prime }(t)$, if $\mathcal{K}=\mathbb{C}^{m}$).

We shall say that a weakly measurable function $t\mapsto \mathbf{C}(t)$ is
locally QS--integrable if its components $C_{\nu }^{\mu },\mu \in
\{-,o\},\;\nu \in \{o,+\}$ are locally $L^{p}$--integrable as
operator--valued functions 
\begin{equation*}
\begin{array}{cc}
C_{+}^{-}(t):\mathcal{G}^{+}\rightarrow \mathcal{G}_{-},\quad & \Vert
C_{+}^{-}(\cdot )\Vert _{\zeta ,t}^{\xi ,1}<\infty \quad (p=1) \\ 
C_{+}^{\circ }(t):\mathcal{G}^{+}\rightarrow \mathcal{G}_{-}\otimes \mathcal{%
E}^{\prime }(t), & \Vert C_{+}^{o}(\cdot )\Vert _{\zeta ,t}^{\xi ,2}<\infty
\quad (p=2) \\ 
C_{\circ }^{-}(t):\mathcal{G}^{+}\otimes \mathcal{E}(t)\rightarrow \mathcal{G%
}_{-}, & \!\Vert C_{\circ }^{o}(\cdot )\Vert _{\zeta ,t}^{\xi ,2}<\infty
\quad (p=2) \\ 
C_{\circ }^{\circ }(t):\mathcal{G}^{+}\otimes \mathcal{E}(t)\rightarrow 
\mathcal{G}_{-}\otimes \mathcal{E}^{\prime }(t),\quad & \ \ \Vert C_{\circ
}^{\circ }(\cdot )\Vert _{\zeta ,t}^{\xi ,\infty }<\infty \quad (p=\infty )%
\end{array}%
\end{equation*}%
Here the norms are defined for any $t>0$, $\xi \in ]0,1[$ and a sufficiently
large $\zeta >1$ by 
\begin{eqnarray*}
\Vert C_{+}^{-}\Vert _{\zeta ,t}^{\xi ,1} &=&\int_{0}^{t}\Vert
C_{+}^{-}(s)\Vert _{\zeta }^{\xi }\mathrm{d}s,\quad \Vert C_{\circ }^{\circ
}\Vert _{\zeta ,t}^{\xi ,\infty }=\mathrm{ess}_{s\leq t}\sup \Vert C_{\circ
}^{\circ }(s)\Vert _{\zeta }^{\xi }, \\
\Vert C_{+}^{-}\Vert _{\zeta }^{\xi } &=&\sup_{\varphi }\{\Vert
C_{+}^{-}\varphi \Vert (\xi )/\Vert \varphi \Vert (\zeta )\},\;\Vert
C_{\circ }^{\circ }\Vert _{\zeta }^{\xi }=\sup_{\varphi ^{\circ }}\{\Vert
C_{\circ }^{\circ }\varphi ^{\circ }\Vert (\xi )/\Vert \varphi ^{\circ
}\Vert (\zeta )\}\ ,
\end{eqnarray*}%
where 
\begin{eqnarray*}
\varphi \in \mathcal{G}(\zeta ),\;\Vert \varphi \Vert ^{2}(\zeta )
&=&<\varphi |\varphi >(\zeta ), \\
\varphi ^{\circ }\in \mathcal{G}(\zeta )\otimes \mathcal{E}(s),\;\Vert
\varphi ^{\circ }\Vert (\zeta ) &=&<\varphi ^{\circ }|\varphi ^{\circ
}>(\zeta )
\end{eqnarray*}%
and $\Vert C_{+}^{\circ }\Vert _{\zeta ,t}^{\xi ,2}=\Vert C_{+}^{\circ
t}\Vert _{\zeta }^{\xi },\;\Vert C_{\circ }^{-}\Vert _{\zeta ,t}^{\xi
,2}=\Vert C_{\circ t}^{-}\Vert _{\zeta }^{\xi }$ are the norms%
\begin{equation*}
\Vert C\Vert _{\xi ,t}^{\xi ,2}=(\int_{0}^{t}(\Vert C(s)\Vert _{\zeta }^{\xi
})^{2}\mathrm{d}s)^{1/2}
\end{equation*}
of the operators 
\begin{eqnarray*}
C_{+}^{\circ t}:\mathcal{G}(\zeta )\rightarrow \mathcal{G}(\xi )\otimes 
\mathcal{E}{^{\prime }}^{t}, &\quad &\left( C_{+}^{\circ t}\varphi \right)
(s)=C_{+}^{\circ }(s)\varphi ,\;s\leq t \\
C_{\circ t}^{-}:\mathcal{G}(\zeta )\otimes \mathcal{E}^{t}\rightarrow 
\mathcal{G}(\xi ), &\quad &C_{\circ t}^{-}\varphi ^{\circ
}=\int_{0}^{t}C_{\circ }^{-}(s)\varphi ^{\circ }(s)\mathrm{d}s
\end{eqnarray*}%
in the Hilbert spaces $\mathcal{E}^{t}=\oplus \int_{0}^{t}\mathcal{E}(s)%
\mathrm{d}s$, $\mathcal{E}{^{\prime }}^{t}=\oplus \int_{0}^{t}\mathcal{E}%
^{\prime }(s)\mathrm{d}s$.

The following theorem shows the continuity of the QS--integral of an
integrable $\mathbf{C}$, defined on $\mathcal{G}^{+}$ even for nonadapted $%
C_{\nu }^{\mu }(t)$ by the formula 
\begin{eqnarray}
\left( \int_{0}^{t}A(\mathbf{C},\mathrm{d}s)\varphi \right) (\tau )
&=&\int_{0}^{t}\left( C_{+}^{-}(s)\varphi +C_{\circ }^{-}(s)\varphi
_{s}^{\circ }\right) (\tau )\mathrm{d}s  \notag \\
&+&\sum_{s\in \tau }^{s\leq t}\left( C_{+}^{\circ }(s)\varphi +C_{\circ
}^{\circ }(s)\varphi _{s}^{\circ }\right) (\tau /s)\ ,  \label{eq:1.10}
\end{eqnarray}%
where $\varphi _{t}^{\circ }\in \mathcal{G}^{+}\otimes \mathcal{E}(t)$ is
defined almost everywhere as the tensor--function $\varphi _{t}^{\circ
}(\tau )=\varphi (\tau \sqcup t)$.

\begin{theorem}
Suppose that $\mathbf{C}(t)$ is a locally QS--integrable function i.e. for
any $\xi <1$, $t>0$ there exists $\zeta >1$, such that 
\begin{equation*}
\Vert C_{+}^{-}\Vert _{\zeta ,t}^{\xi ,1}<\infty ,\;\Vert C_{+}^{\circ
}\Vert _{\zeta ,t}^{\xi ,2}<\infty ,\;\Vert C_{\circ }^{-}\Vert _{\zeta
,t}^{\xi ,2}<\infty ,\;\Vert C_{\circ }^{\circ }\Vert _{\zeta ,t}^{\xi
,\infty }<\infty .
\end{equation*}%
Then the QS--integral (\ref{eq:1.10}) is defined as a continuous operator $%
\imath _{0}^{t}(\mathbf{C}):\mathcal{G}^{+}\rightarrow \mathcal{G}_{-}$ with
the estimate 
\begin{equation}
\left\Vert \int_{0}^{t}A(\mathbf{C}(s),\mathrm{d}s)\right\Vert _{\eta
^{+}}^{\eta _{-}}\leq \Vert C_{+}^{-}\Vert _{\zeta ,t}^{\xi ,1}+{\frac{1}{%
\sqrt{\varepsilon }}}\left( \Vert C_{\circ }^{-}\Vert _{\zeta ,t}^{\xi
,2}+\Vert C_{+}^{\circ }\Vert _{\zeta ,t}^{\xi ,2}\right) +{\frac{1}{%
\varepsilon }}\Vert C_{\circ }^{\circ }\Vert _{\zeta ,t}^{\xi ,\infty }
\label{eq:1.11}
\end{equation}%
for the norms $\Vert \imath _{0}^{t}(\mathbf{C})\Vert _{\eta ^{+}}^{\eta
_{-}}=\sup_{\varphi }\{\Vert \imath _{0}^{t}(C)\varphi \Vert (\eta
_{-})/\Vert \varphi \Vert (\eta ^{+})\}$, where $\eta _{-}\leq \xi
-\varepsilon $, $\eta ^{+}\geq \zeta +\varepsilon $ and $0<\varepsilon <\xi $%
. Moreover, the adjoint integral 
\begin{equation*}
<\int_{0}^{t}A(\mathbf{C},\mathrm{d}s)^{\ast }\varphi |\chi >=<\varphi
|\int_{0}^{t}A(\mathbf{C},\mathrm{d}s)\chi >,\;\varphi ,\chi \in \mathcal{G}%
^{+}
\end{equation*}%
is also densely defined on $\mathcal{G}^{+}\subset \mathcal{G}_{-}$ as the
QS--integral $\int_{0}^{t}A(\mathbf{C}^{\star },\mathrm{d}s)$, and the
function $\mathbf{C}^{\star }(t)=\mathbf{g}\mathbf{C}(t)^{\ast }\mathbf{g}$, 
\begin{equation*}
\mathbf{C}^{\star }(t)_{+}^{-}=C_{+}^{-}(t)^{\ast },\;\mathbf{C}^{\star
}(t)_{+}^{\circ }=C_{\circ }^{-}(t)^{\ast },\;\mathbf{C}^{\star }(t)_{\circ
}^{-}=C_{+}^{\circ }(t)^{\ast },\;\mathbf{C}^{\star }(t)_{\circ }^{\circ
}=C_{\circ }^{\circ }(t)^{\ast }
\end{equation*}%
is locally QS--integrable with $\Vert \mathbf{C}^{\star }(t)_{\nu }^{\mu
}\Vert _{1/\xi }^{1/\zeta }=\Vert C_{\nu }^{\mu }(t)^{\ast }\Vert _{\zeta
}^{\xi }<\infty $ for almost all $t$. \label{th:1}
\end{theorem}

\noindent \textsc{Proof.}\ . In order to show the continuity of the integral
(\ref{eq:1.10}) in the projective topology of $\bigcap_{\zeta >1}\mathcal{G}%
(\zeta )$, one should prove that 
\begin{equation*}
\Vert \int_{0}^{t}A(\mathbf{C}(s),\mathrm{d}s)\varphi \Vert (\eta _{-})\leq
c\Vert \varphi \Vert (\eta ^{+}),\Vert \varphi \Vert (\eta )=(<\varphi
|\varphi >(\eta ))^{1/2}
\end{equation*}%
for any $\varphi \in \mathcal{G}(\eta ^{+})$, $\eta _{-}<\xi $ and a $\eta
^{+}>\zeta ,\;c>0$. Due to the definition 
\begin{eqnarray*}
\left\Vert \int_{0}^{t}A(\mathbf{C},\mathrm{d}s)\varphi \right\Vert &\leq
&\Vert \int_{0}^{t}C_{\circ }^{-}\mathrm{d}A_{-}^{\circ }\varphi \Vert
+\Vert \int_{0}^{t}C_{+}^{\circ }\mathrm{d}A_{\circ }^{+}\varphi \Vert \\
&&+\Vert \int_{0}^{t}C_{\circ }^{\circ }\mathrm{d}N_{\circ }^{\circ }\varphi
\Vert +\Vert \int_{0}^{t}C_{+}^{-}\mathrm{d}s\varphi \Vert ,
\end{eqnarray*}%
where 
\begin{equation*}
\int_{0}^{t}C_{\circ }^{-}\mathrm{d}A_{-}^{\circ }\varphi
=\int_{0}^{t}C_{\circ }^{-}(s)\varphi _{s}^{\circ }\mathrm{d}s,\;\left(
\int_{0}^{t}C_{+}^{\circ }\mathrm{d}A_{\circ }^{+}\varphi \right) (\tau
)=\sum_{s\in \tau }^{s\leq t}(C_{+}^{\circ }(s)\varphi )(\tau \backslash s),
\end{equation*}%
and 
\begin{equation*}
\left( \int_{0}^{t}C_{\circ }^{\circ }\mathrm{d}N_{\circ }^{\circ }\varphi
\right) (\tau )=\sum_{s\in \tau }^{s\leq t}\left( C_{\circ }^{\circ
}(s)\varphi _{s}^{\circ }\right) (\tau \backslash s).
\end{equation*}%
The first two integrals in (\ref{eq:1.10}) can be easily estimated as 
\begin{eqnarray*}
\biggl\Vert\int_{0}^{t}C_{+}^{-}\varphi \mathrm{d}s\biggr\Vert(\eta _{-})
&\leq &\int_{0}^{t}\Vert C_{+}^{-}(s)\varphi \Vert (\xi )\mathrm{d}s\leq
\int_{0}^{t}\Vert C_{+}^{-}(s)\Vert _{\zeta }^{\xi }\mathrm{d}s\Vert \varphi
\Vert (\zeta )= \\
&&\Vert C_{+}^{-}\Vert _{\zeta ,t}^{\xi ,1}\Vert \varphi \Vert (\zeta ), \\
\biggl\Vert\int_{0}^{t}C_{\circ }^{-}\mathrm{d}A_{-}^{\circ }\varphi %
\biggr\Vert(\eta _{-}) &=&\Vert C_{\circ t}^{-}\varphi ^{\circ }\Vert (\xi
)\leq \Vert C_{\circ t}^{-}\Vert _{\zeta }^{\xi }\Vert \varphi ^{\circ
}\Vert (\zeta )= \\
&&\Vert C_{\circ }^{-}\Vert _{\zeta ,t}^{\xi ,2}\left( {\frac{\mathrm{d}}{%
\mathrm{d}\zeta }}\Vert \varphi \Vert ^{2}(\zeta )\right) ^{1/2}\ ,
\end{eqnarray*}%
where we took into account that 
\begin{equation*}
\Vert \varphi ^{\circ }\Vert ^{2}(\zeta )=\int \int \zeta ^{|\tau |}\Vert
\varphi (\tau \sqcup t)\Vert ^{2}\mathrm{d}\tau \mathrm{d}t=\int |\tau
|\zeta ^{|\tau |-1}\Vert \varphi (\tau )\Vert ^{2}\mathrm{d}\tau ={\frac{%
\mathrm{d}}{\mathrm{d}\zeta }}\Vert \varphi \Vert ^{2}(\zeta ).
\end{equation*}%
In order to estimate the integrals of $C_{+}^{\circ }$ and $C_{\circ
}^{\circ }$ let us find 
\begin{eqnarray*}
\lefteqn{\biggl\Vert \int_{0}^{t}C_{+}^{\circ }\mathrm{d}A_{\circ
}^{+}\varphi \biggr\Vert ^{2}(\eta _{-})=\int_{\Omega }\left\Vert \sum_{s\in
\tau }^{s\leq t}(C_{+}^{\circ }(s)\varphi )(\tau /s)\right\Vert ^{2}\eta
_{-}^{|\tau |}\mathrm{d}\tau =} \\
&&\eta _{-}^{2}\int_{0}^{t}\int_{0}^{t}\int_{\Omega }<[C_{+}^{\circ
}(s_{1})\varphi _{s_{2}}^{\circ }](\tau )|[C_{+}^{\circ }(s_{2})\varphi
_{s_{1}}^{\circ }](\tau )>\eta _{-}^{|\tau |}\mathrm{d}\tau \mathrm{d}s_{1}%
\mathrm{d}s_{2}+ \\
&&\eta _{-}\int_{0}^{t}\int_{\Omega }\Vert \lbrack C_{+}^{\circ }(s)\varphi
](\tau )\Vert ^{2}\eta _{-}^{|\tau |}\mathrm{d}\tau \mathrm{d}s\leq \eta
_{-}\left( 1+\eta _{-}{\frac{\mathrm{d}}{\mathrm{d}\eta _{-}}}\right) \Vert
C_{+}^{\circ }\varphi \Vert _{t}^{2}(\eta _{-})
\end{eqnarray*}%
by Schwarz inequality. In the same way we get 
\begin{eqnarray*}
\lefteqn{\biggl\Vert \int_{0}^{t}C_{\circ }^{\circ }\mathrm{d}N_{\circ
}^{\circ }\varphi \biggr\Vert ^{2}(\eta _{-})=\int_{\Omega }\left\Vert
\sum_{s\in \tau }^{s\leq t}(C_{\circ }^{\circ }(s)\varphi ^{\circ }(s))(\tau
/s)\right\Vert ^{2}\eta _{-}^{|\tau |}\mathrm{d}\tau =} \\
&=&\eta _{-}^{2}\int_{0}^{t}\int_{0}^{t}\int_{\Omega }<[C_{\circ }^{\circ
}(s_{1})\varphi _{s_{1}\,s_{2}}^{\circ \;\,\circ }](\tau )|[C_{\circ
}^{\circ }(s_{2})\varphi _{s_{1}\,s_{2}}^{\circ \;\,\circ }](\tau )>\eta
_{-}^{|\tau |}\mathrm{d}\tau \mathrm{d}s_{1}\mathrm{d}s_{2}+ \\
&&\eta _{-}\int_{0}^{t}\int_{\Omega }\Vert \lbrack C_{\circ }^{\circ
}(s)\varphi _{s}^{\circ }](\tau )\Vert ^{2}\eta _{-}^{|\tau |}\mathrm{d}\tau 
\mathrm{d}s\leq \eta _{-}\left( 1+\eta _{-}{\frac{\mathrm{d}}{\mathrm{d}\eta
_{-}}}\right) \Vert C_{\circ }^{\circ }\varphi ^{\circ }\Vert _{t}^{2}(\eta
_{-})\ ,
\end{eqnarray*}%
where $\varphi _{s_{1},s_{2}}^{\circ \circ }(\tau )=\varphi (\tau \sqcup
s_{1}\sqcup s_{2})$.

Taking into account that for any $\varepsilon >0$, $\xi =\eta +\varepsilon $ 
\begin{equation*}
{\frac{\mathrm{d}}{\mathrm{d}\zeta }}\Vert \varphi \Vert ^{2}(\eta )\leq {%
\frac{1}{\varepsilon }}\left( \Vert \varphi \Vert ^{2}(\eta +\varepsilon
)-\Vert \varphi \Vert ^{2}(\eta )\right) \leq {\frac{1}{\varepsilon }}\Vert
\varphi \Vert ^{2}(\xi )\ ,
\end{equation*}%
one can find that $\left( 1+\eta \,{\frac{\mathrm{d}}{\mathrm{d}\eta }}%
\right) \Vert \varphi \Vert _{t}^{2}(\eta )\leq {\frac{\xi }{\varepsilon }}%
\Vert \varphi \Vert _{t}^{2}(\xi )$, 
\begin{eqnarray*}
\left( 1+\eta {\frac{\mathrm{d}}{\mathrm{d}\eta }}\right) \Vert C_{+}^{\circ
}\varphi \Vert _{t}^{2}(\eta ) &\leq &{\frac{\xi }{\varepsilon }}\left(
\Vert C_{+}^{\circ }\Vert _{\zeta ,t}^{\xi ,2}\right) ^{2}\Vert \varphi
\Vert ^{2}(\zeta )\ , \\
\left( 1+\eta {\frac{\mathrm{d}}{\mathrm{d}\eta }}\right) \Vert C_{\circ
}^{\circ }\varphi ^{\circ }\Vert _{t}^{2}(\eta ) &\leq &{\frac{\xi }{%
\varepsilon }}\left( \Vert C_{\circ }^{\circ }\Vert _{\zeta ,t}^{\xi ,\infty
}\right) ^{2}{\frac{\mathrm{d}}{\mathrm{d}\zeta }}\Vert \varphi \Vert
^{2}(\zeta )\ ,
\end{eqnarray*}%
if $\varepsilon \leq \xi $. Hence, due to $\Vert \varphi (\eta ^{+})\geq
\Vert \varphi \Vert (\zeta +\varepsilon )\geq \Vert \varphi \Vert (\zeta )$
for $\eta ^{+}\geq \zeta +\varepsilon $, we obtain for $\eta _{-}\leq \eta
=\xi -\varepsilon $, $\xi \leq 1$ 
\begin{equation*}
\left\Vert \int_{0}^{t}A(\mathbf{C},\mathrm{d}s)\varphi \right\Vert (\eta
^{-})\leq \Vert C_{+}^{-}\Vert _{\zeta ,t}^{\xi ,1}+{\frac{1}{\sqrt{%
\varepsilon }}}\left( \Vert C_{\circ }^{-}\Vert _{\zeta ,t}^{\xi ,2}+\Vert
C_{+}^{\circ }\Vert _{\zeta ,t}^{\xi ,2}\right) +{\frac{1}{\varepsilon }}%
\Vert C_{\circ }^{\circ }\Vert _{\zeta ,t}^{\xi ,\infty },
\end{equation*}%
if $\Vert \varphi \Vert (\eta ^{+})\leq 1,\;0<\varepsilon \leq \xi $, what
is equivalent to (\ref{eq:1.11}).

Due to the duality $\mathcal{G}(\zeta )^{\ast }=\mathcal{G}(\zeta ^{-1})$ of 
$\mathcal{F}(\zeta )$ and $\mathcal{G}(1/\zeta )$ the QS--matrix process $%
\mathbf{C}^{\star }(t)$ is also locally QS--integrable, and there exists the
adjoint integral $\int_{0}^{t}A(\mathbf{C},\mathrm{d}s)^{\ast }$, defined as
in (\ref{eq:1.10}) by $\mathbf{C}^{\star }$: 
\begin{eqnarray*}
\lefteqn{<\varphi |\int_{0}^{t}A(C,\mathrm{d}s)\chi >=\int_{0}^{t}<\varphi
|C_{+}^{-}(s)\chi +C_{\circ }^{-}(s)\chi ^{\circ }(s)>\mathrm{d}s} \\
&+&\int_{\circ }^{t}<\varphi _{s}^{\circ }|C_{+}^{\circ }(s)\chi +C_{\circ
}^{\circ }(s)\chi _{s}^{\circ }>\mathrm{d}s=\int_{0}^{t}<C_{+}^{-}(s)^{\ast
}\varphi +C_{+}^{\circ }(s)^{\ast }\varphi _{s}^{\circ }|\chi >\mathrm{d}s+
\\
&+&\int_{0}^{t}<C_{\circ }^{-}(s)^{\ast }\varphi +C_{\circ }^{\circ
}(s)^{\ast }\varphi _{s}^{\circ }|\chi ^{\circ }>\mathrm{d}s=<\int_{0}^{t}A(%
\mathbf{C}^{\star },\mathrm{d}s)\varphi |\chi >\ .
\end{eqnarray*}%
Obviously, $\Vert \int_{0}^{+}A(\mathbf{C}^{\star },\mathrm{d}s)\Vert
_{1/\eta _{-}}^{1/\eta _{+}}=\Vert \int_{0}^{+}A(\mathbf{C},\mathrm{d}%
s)\Vert _{\eta _{+}}^{\eta _{-}}$\hfill \vrule height.9exwidth.8exdepth-.1ex

\begin{corollary}
If $\mathbf{C}(t)$ are the simple measurable adapted functions, then the
definition (\ref{eq:1.10}) coincides with the QS--integral, given by
integral Ito's sums with respect to the processes (\ref{eq:1.1}). Moreover,
the QS--integral (\ref{eq:1.10}) is a limit of such integral sums in the
inductive operator topology, defined by the norms (\ref{eq:1.11}), if
locally QS--integrable matrix--process $\mathbf{C}$ can be uniformly
approximated by a sequence of simple operator--valued processes with respect
to the defined $L^{p}$ -- norms on $]0,t]$. \label{corol:1}
\end{corollary}


\section{QS calculus of output nondemolition processes}

\label{sec:belat2}

Let us consider an initial Hilbert space $\mathcal{H}^{0}=\mathfrak{h}$ with
identity operator $\widehat{1}$, $\mathcal{H}=\mathfrak{h}\otimes \mathcal{G}
$, and denote by $\mathcal{H}^{t}=\mathfrak{h}\otimes \mathcal{G}^{t}$ and
by $\widehat{I}^{t}=\widehat{1}\otimes I^{t}$ the corresponding multipliers
of the Hilbert space $\mathcal{H}=\mathcal{H}^{t}\otimes \mathcal{G}_{t}$
and identity operator $\widehat{I}=\widehat{I}^{t}\otimes \widehat{1}_{t}$.
Let us identify the basic QS process $\mathbf{A}=\left( A_{\nu }^{\mu
}\right) $ with the process $\widehat{\mathbf{A}}=\widehat{1}\otimes \mathbf{%
A}$, in $\mathcal{H}$: $\widehat{A}_{+}^{-}(t)=t\widehat{I},\;\widehat{A}%
_{-}^{j}=\widehat{1}\otimes A_{-}^{j},\;\widehat{A}_{j}^{+}=\widehat{1}%
\otimes A_{j}^{+},\;\widehat{A}_{k}^{i}=\widehat{1}\otimes N_{k}^{i}$. A QS
matrix process $\widehat{\mathbf{C}}=\left( \widehat{C}_{\nu }^{\mu }\right) 
$ with an $\widehat{C}_{\nu }^{\mu }(t)$ acting in $\mathcal{H}$ is called
adapted, if $\widehat{\mathbf{C}}(t)=\widehat{\mathbf{C}}^{t}\otimes 
\widehat{I}_{t}$, for any $t$, where $\widehat{\mathbf{C}}^{t}$ is a matrix
of operators in $\mathcal{H}^{t}$. We define the QS integral of an adapted
QS matrix process $\widehat{\mathbf{C}}$ as in (\ref{eq:1.10}) by the sum of
integrals 
\begin{equation}
\int_{0}^{t}\widehat{A}\left( \widehat{\mathbf{C}},\mathrm{d}s\right)
=\int_{0}^{t}\left( \widehat{C}_{+}^{-}\mathrm{d}s+\widehat{C}_{\circ }^{-}%
\mathrm{d}\widehat{A}_{-}^{\circ }+\widehat{C}_{-}^{\circ }\mathrm{d}%
\widehat{A}_{\circ }^{+}+\widehat{C}_{\circ }^{\circ }\mathrm{d}\widehat{N}%
_{\circ }^{\circ }\right) \ ,  \label{eq:2.1}
\end{equation}%
which exists as an adapted process with the QS differential $\widehat{A}(%
\widehat{\mathbf{C}},\mathrm{d}t)=\widehat{C}_{\nu }^{\mu }(t)\mathrm{d}%
\widehat{A}_{\mu }^{\nu }(t)$ for weakly measurable, locally integrable
functions $t\mapsto \widehat{C}_{\nu }^{\mu }(t)$, called below QS
integrable processes.

Now let us consider an adapted process $\widehat{X}(t)$, defined by the QS
differential equation 
\begin{equation*}
\mathrm{d}\widehat{X}(t)=\left( \widehat{F}_{\nu }^{\mu }(t)-\widehat{X}%
(t)\delta _{\nu }^{\mu }\right) \mathrm{d}\widehat{A}_{\mu }^{\nu }(t)\
,\qquad \qquad \widehat{X}(0)=\widehat{x}\otimes I,
\end{equation*}%
having the solution $\widehat{X}(t)=\widehat{X}\otimes I+\int_{0}^{t}%
\widehat{C}_{\nu }^{\mu }\mathrm{d}\widehat{A}_{\mu }^{\nu }$, iff $\widehat{%
\mathbf{C}}=\widehat{\mathbf{F}}-\widehat{X}\otimes \mathbf{\delta }$
satisfies the conditions for the existence of the integral (\ref{eq:2.1}),
where $\widehat{X}\otimes \mathbf{\delta }=\left( \widehat{X}\delta _{\nu
}^{\mu }\right) $. We shall define the elements $\widehat{F}_{\nu }^{\mu }$
of matrix--operators $\widehat{\mathbf{F}}(t)$ also for $\mu =-=\nu $ and
for $\mu =+=\nu $ by $\widehat{F}_{-}^{-}=\widehat{X}=\widehat{F}_{+}^{+}$,
and assume that $\widehat{F}_{\nu }^{\mu }=0$, if $\mu >\nu $ under the
order $-<o<+$.

\begin{proposition}
If the QS process $\widehat{X}$ satisfies the QS differential equation (\ref%
{eq:2.2}), then the process $(\widehat{X}^{\ast }\widehat{X})(t)=\widehat{X}%
(t)^{\ast }\widehat{X}(t)$ satisfies the equation 
\begin{equation}
\mathrm{d}(\widehat{X}^{\ast }\widehat{X})=\left( \widehat{\mathbf{F}}%
^{\star }\widehat{\mathbf{F}}-\widehat{X}^{\ast }\widehat{X}\otimes \mathbf{%
\delta }\right) _{\nu }^{\mu }\mathrm{d}A_{\mu }^{\nu },\qquad \qquad (%
\widehat{X}^{\ast }\widehat{X})(0)=\widehat{x}^{\ast }\widehat{x}\otimes I.
\label{eq:2.3}
\end{equation}%
This QS Ito formula establishes an $\ast $-algebra isomorphism from the QS
differentiable processes $\widehat{X}$ into the algebra of matrices of
operator processes $\widehat{F}_{\nu }^{\mu }$ defined above. In particular, 
$\widehat{X}$ is formally normal (selfadjoint, unitary) iff $[\widehat{%
\mathbf{F}},\widehat{\mathbf{F}}^{\star }]=0,\left( \widehat{\mathbf{F}}%
^{\star }=\widehat{\mathbf{F}},\widehat{\mathbf{F}}^{\star }=\widehat{%
\mathbf{F}}^{-1}\right) $ with repsect to the $\star $-operation $\mathbf{F}%
^{\star }(t)=\mathbf{g}\mathbf{F}(t)^{\ast }\mathbf{g}$, and $X$ is
partially isometric (isometric, orthoprojection), iff $\widehat{\mathbf{F}}%
\widehat{\mathbf{F}}^{\star }\widehat{\mathbf{F}}=\widehat{\mathbf{F}}$ ($%
\widehat{\mathbf{F}}^{\star }\widehat{\mathbf{F}}=\widehat{I}\otimes \mathbf{%
\delta },\widehat{\mathbf{F}}^{\star }\widehat{\mathbf{F}}=\widehat{\mathbf{F%
}}$). \label{prop:2}
\end{proposition}

\noindent \textsc{Proof.}\ Taking into account that 
\begin{eqnarray*}
\mathrm{d}\widehat{X} &=&\widehat{C}_{+}^{-}\mathrm{d}t+\widehat{C}_{\circ
}^{-}\mathrm{d}\widehat{A}_{-}^{\circ }+\widehat{C}_{+}^{\circ }\mathrm{d}%
\widehat{A}_{\circ }^{+}+\widehat{C}_{\circ }^{\circ }\mathrm{d}\widehat{N}%
_{\circ }^{\circ }, \\
\mathrm{d}\widehat{X}^{\ast } &=&\widehat{C}_{+}^{-\ast }\mathrm{d}t+%
\widehat{C}_{+}^{\circ \ast }\mathrm{d}\widehat{A}_{-}^{\circ }+\widehat{C}%
_{\circ }^{-\ast }\mathrm{d}\widehat{A}_{\circ }^{+}+\widehat{C}_{\circ
}^{\circ \ast }\mathrm{d}\widehat{N}_{\circ }^{\circ },
\end{eqnarray*}%
and using the QS It\^{o} formula \cite{bib:12}, defining the product $(%
\widehat{X}^{\ast }\widehat{X})(t)=(\widehat{x}^{\ast }\widehat{x})\otimes
I+\int_{0}^{t}\mathrm{d}(\widehat{X}^{\ast }\widehat{X})$ by the QS
differential 
\begin{equation*}
\mathrm{d}(\widehat{X}^{\ast }\widehat{X})=\mathrm{d}\widehat{X}^{\ast }%
\widehat{X}+\widehat{X}^{\ast }\mathrm{d}\widehat{X}+\mathrm{d}\widehat{X}%
^{\ast }\mathrm{d}\widehat{X}=
\end{equation*}%
\begin{equation*}
=\left( C_{\;\;\nu }^{\star \mu }\widehat{X}+\widehat{X}^{\ast }\widehat{C}%
_{\nu }^{\mu }+C_{\;\;\mu }^{\star \mu }C_{\nu }^{\mu }\right) \mathrm{d}%
\widehat{A}_{\mu }^{\nu }=
\end{equation*}%
\begin{equation*}
\left( (\widehat{\mathbf{C}}+\widehat{X}\otimes \mathbf{\delta })^{\star }(%
\widehat{\mathbf{C}}+\widehat{X}\otimes \mathbf{\delta })-\widehat{X}^{\ast }%
\widehat{X}\otimes \mathbf{\delta }\right) _{\nu }^{\mu }\mathrm{d}\widehat{A%
}_{\mu }^{\nu },
\end{equation*}%
we obtain the equation (\ref{eq:2.3}) with $\widehat{\mathbf{F}}=\widehat{X}%
\otimes \mathbf{\delta }+\widehat{\mathbf{C}}$. Due to the linearity of (\ref%
{eq:2.3}) with respect to the pairs $(\widehat{\mathbf{F}},\widehat{X})$ and 
$(\widehat{\mathbf{F}}^{\star },\widehat{X}^{\ast })$, it can be extended to 
\begin{equation}
\mathrm{d}(\widehat{X}^{\ast }\widehat{X}^{\prime })=\left( \widehat{\mathbf{%
F}}^{\star }\widehat{\mathbf{F}}^{\prime }-\widehat{X}^{\ast }\widehat{X}%
^{\prime }\otimes \mathbf{\delta }\right) _{\nu }^{\mu }\mathrm{d}\widehat{A}%
_{\mu }^{\nu }  \label{eq:2.4}
\end{equation}%
by the polarization formula 
\begin{equation*}
\widehat{X}^{\ast }\widehat{X}^{\prime }=\sum_{n=0}^{3}\left( \widehat{X}+%
\mathrm{i}^{n}\widehat{X}^{\prime }\right) ^{\ast }\left( \widehat{X}+%
\mathrm{i}^{n}\widehat{X}^{\prime }\right) /4\mathrm{i}^{n}\ ,\qquad \qquad 
\mathrm{i}=\sqrt{-1}.
\end{equation*}%
Hence, the formula (\ref{eq:2.3}) is equivalent to QS Hudson --
Parthasarathy It\^{o} formula \cite{bib:12} and $\star $ -- property 
\begin{equation*}
\mathrm{d}\widehat{X}^{\ast }(t)=\left( \widehat{\mathbf{F}}^{\star }(t)-%
\widehat{X}^{\ast }(t)\otimes \mathbf{\delta }\right) _{\nu }^{\mu }\mathrm{d%
}\widehat{A}_{\mu }^{\nu },\qquad \qquad \widehat{X}^{\ast }(0)=\widehat{x}%
^{\ast }\otimes I,
\end{equation*}%
%
%
%
%
%
%
%
%
which follows from it for $\widehat{\mathbf{F}}^{\prime }=I\otimes \mathbf{%
\delta }$, corresponding to $\widehat{X}^{\prime }=\widehat{I}$. So the map $%
\widehat{X}\mapsto \widehat{\mathbf{F}}$ is a homomorphism with respect to
the associative operator algebra structure of $\widehat{X}$ and $\widehat{%
\mathbf{F}}$ with the appropriate involutions. Furthermore it is an
injection, because if $\widehat{\mathbf{F}}=0$, then $\widehat{X}=0$, as $%
\widehat{F}_{-}^{-}=\widehat{X}=\widehat{F}_{+}^{+}$.

Conversely, if $\widehat{X}=0$, then $\int_{0}^{t}\widehat{F}_{\nu }^{\mu }%
\mathrm{d}\widehat{A}_{\mu }^{\nu }=0$ for all $t$, but it implies $\widehat{%
F}_{\nu }^{\mu }=0$ due to the independence of stochastic integrators \cite%
{bib:belat4}.\hfill \vrule height .9ex width .8ex depth -.1ex

Now let us consider an adapted selfadjoint QS process $Y$, satisfying a QS
equation 
\begin{equation}
\mathrm{d}Y(t)=\left( \mathbf{Z}^{\star }\mathbf{GZ}-Y\otimes \mathbf{\delta 
}\right) _{\nu }^{\mu }(t)\mathrm{d}\widehat{A}_{\mu }^{\nu }(t)\ ,\qquad
\qquad Y(0)=\widehat{y}\otimes I\ ,  \label{eq:2.5}
\end{equation}%
where $\mathbf{G}^{\star }=\mathbf{G}$ is a $\star $ -- selfadjoint matrix
adapted QS process with $G_{-}^{-}=Y=G_{+}^{+},\,G_{\nu }^{\mu }=0$ for $\mu
>\nu $, and $\mathbf{Z}=\left( Z_{\nu }^{\mu }\right) $ is a $\star $ --
isometric or $\star $ -- unitary matrix adapted process: $\mathbf{Z}^{\star }%
\mathbf{Z}=\widehat{I}\otimes \mathbf{\delta }\,\left( \mathbf{Z}^{\star }=%
\mathbf{Z}^{-1}\right) $ such, that $\mathbf{GZZ}^{\star }=\mathbf{G}=%
\mathbf{ZZ}^{\star }\mathbf{G}$ (otherwise $\mathbf{G}$ should be replaced
by $\mathbf{ZZ}^{\star }\mathbf{GZZ}^{\star }$).

We shall demand that $Z_{-}^{-}=\widehat{I}=Z_{+}^{+},\,Z_{\nu }^{\mu }=0$,
if $\mu >\nu $, and $Z_{\nu }^{\mu },\mu \not=+$ or $\nu \not=-$ satisfy the
conditions for the existance of QS isometric (unitary) evolution $U(t):%
\mathcal{H}\rightarrow \mathcal{H}$, defined by the QS equation (\ref{eq:2.2}%
) with $\widehat{X}=U,F_{\nu }^{\mu }=UZ_{\nu }^{\mu },\widehat{x}=\widehat{1%
}$: 
\begin{equation}
\mathrm{d}U(t)=U(t)\left( Z_{\nu }^{\mu }(t)-I\delta _{\nu }^{\mu }\right) 
\mathrm{d}\widehat{A}_{\mu }^{\nu }\ ,\qquad \qquad U(0)=\widehat{I}.
\label{eq:2.6}
\end{equation}%
Sufficient conditions for this are the conditions of local integrability of
the weakly measurable processes $Z_{\nu }^{\mu }$ in the sense of the $L^{p}$%
-norms \cite{bib:belat5}: 
\begin{equation*}
\Vert Z_{+}^{-}\Vert _{t}^{(1)}<\infty \ ,\;\Vert Z_{+}^{0}\Vert
_{t}^{(2)}<\infty \ ,\;\Vert Z_{0}^{-}\Vert _{t}^{(2)}<\infty \ ,\;\Vert
Z_{0}^{0}\Vert _{t}^{(\infty )}<\infty \ .
\end{equation*}

Let us call the process $Y$ an output process, if $Y$ is nondemolition with
respect to the QS process $\mathbf{Z}$, generating the evolution (\ref%
{eq:2.6}). By this we mean the commutativity condition 
\begin{equation}
\lbrack Y(t),X(s)]=0,\qquad \forall t\leq s  \label{eq:2.7}
\end{equation}%
with respect to the all QS processes $X(t)=Z_{\nu }^{\mu }(t),\mu ,\nu \in
\{-,J,+\}$ (the conditions are nontrivial for $\mu \not=+$ and $\nu \not=-$).

\begin{theorem}
The process $Y$ is defined by (\ref{eq:2.5}) as an adapted selfadjoint QS
process iff 
\begin{equation}
UY=\widehat{Y}U,\quad \widehat{Y}(t)=\widehat{y}\otimes I+\int_{0}^{t}%
\widehat{D}_{\nu }^{\mu }\mathrm{d}\widehat{A}_{\mu }^{\nu },  \label{eq:2.8}
\end{equation}%
where $\widehat{\mathbf{D}}$ is an adapted QS integrable matrix process
satisfying the conditions $\widehat{D}_{\nu }^{\mu }U=UD_{\nu }^{\mu
},D_{\nu }^{\mu }=G_{\nu }^{\mu }-Y\delta _{\nu }^{\mu }$. The process $Y$
is an output QS process, iff 
\begin{equation}
U(s)Y(t)=\widehat{Y}(t)U(s)\quad ,\quad \forall t\leq s,  \label{eq:2.9}
\end{equation}%
which is equivalent to the condition $[\widehat{Y}(t),\widehat{Z}_{\nu
}^{\mu }(s)]U(s)=0$ for $s\leq t,\widehat{Z}_{\nu }^{\mu }U=UZ_{\nu }^{\mu }$%
. The output process $Y$ satisfies the nondemolition condition (\ref{eq:2.7}%
) with respect to an adapted QS process $X$, defined by 
\begin{equation}
UX=\widehat{X}U,\quad \widehat{X}(t)=\widehat{x}\otimes I+\int_{0}^{t}%
\widehat{C}_{\nu }^{\mu }\mathrm{d}\widehat{A}_{\mu }^{\nu },
\label{eq:2.10}
\end{equation}%
iff $[\widehat{Y}(t),\widehat{X}(s)]U(s)=0$ for all $t\leq s$. The last is
equivalent to the commutativity conditions $[\widehat{Y}(t),\widehat{F}_{\nu
}^{\mu }(s)]U(s)=0,\forall t\leq s$, 
\begin{equation}
\lbrack \widehat{y},\widehat{x}]=0,\quad \lbrack \widehat{\mathbf{D}},%
\widehat{\mathbf{F}}]U=0,  \label{eq:2.11}
\end{equation}%
where $\widehat{\mathbf{F}}=\widehat{\mathbf{C}}+\widehat{X}\otimes \mathbf{%
\delta }$, (the conditions are nontrivial for $\mu \not=+$ and $\nu \not=-$).%
\label{th:2}
\end{theorem}

\noindent \textsc{Proof.}\ We obtain (\ref{eq:2.5}) with $G_{\nu }^{\mu
}=U^{\ast }\widehat{G}_{\nu }^{\mu }U$ from (\ref{eq:2.8}) for $U$,
satisfying (\ref{eq:2.6}) simply by applying to $Y=U^{\ast }\widehat{Y}U$
the QS It\^{o} formula (\ref{eq:2.5}): 
\begin{equation*}
\mathrm{d}(U^{\ast }\widehat{Y}U)=\left( \mathbf{Z}^{\star }(U^{\ast
}\otimes \mathbf{\delta })\widehat{\mathbf{G}}(U\otimes \mathbf{\delta })%
\mathbf{Z}-U^{\ast }\widehat{Y}U\otimes \mathbf{\delta }\right) _{\nu }^{\mu
}\mathrm{d}\widehat{A}_{\mu }^{\nu }.
\end{equation*}%
Conversely, we obtain 
\begin{equation*}
\mathrm{d}(UYU^{\ast })=\left( (U\otimes \mathbf{\delta })\mathbf{ZZ}^{\star
}\mathbf{GZZ}^{\star }(U^{\ast }\otimes \mathbf{\delta })-UYU^{\ast }\otimes 
\mathbf{\delta }\right) _{\nu }^{\mu }\mathrm{d}\widehat{A}_{\mu }^{\nu },
\end{equation*}%
so the QS--process $\widehat{Y}=UYU^{\ast }$ obviously satisfying the
condition $\widehat{Y}U=UY$, is defined as QS integral in (\ref{eq:2.8})
with $G_{\nu }^{\mu }U^{\ast }$ due to the assumption $\mathbf{ZZ}^{\star }%
\mathbf{GZZ}^{\star }=\mathbf{G}$ for a weakly measurable, locally $L^{p}$%
-integrable $\mathbf{G}$.

If the processes $Y$ and $X=Z_{\nu }^{\mu }$ satisfy the commutativity
conditions (\ref{eq:2.7}), then the isometry 
\begin{equation}
U(t,s)=\widehat{I}+\int_{t}^{s}U(t,r)\left( Z_{\nu }^{\mu }(r)-\widehat{I}%
\delta _{\nu }^{\mu }\right) \mathrm{d}\widehat{A}_{\mu }^{\nu }(r)
\label{eq:2.12}
\end{equation}%
commutes with $Y(t)$, as it can be easily proved by induction with respect
to $n=1,2,\dots $ for the corresponding QS It\^{o} integral sums 
\begin{equation}
U_{n}(t,s)=\widehat{I}+\sum_{i=0}^{n-1}U_{i}(t,t_{i})\left( Z_{\nu }^{\mu
}(t_{i})-\widehat{I}\delta _{\nu }^{\mu }\right) \left( \widehat{A}_{\mu
}^{\nu }\left( t_{i+1}\right) -\widehat{A}_{\mu }^{\nu }(t_{i})\right) ,
\label{eq:2.13}
\end{equation}%
where $t_{i}=t+i(s-t)/n,U_{0}(t,t)=\widehat{I}$. Hence, taking into account
that $U(s)=U(t)U(t,s)$, we obtain (\ref{eq:2.5}): 
\begin{equation*}
U(s)Y(t)=U(t)Y(t)U(t,s)=\widehat{Y}(t)U(s).
\end{equation*}%
Conversely, multiplying (\ref{eq:2.9}) from the left side hand by $%
U(t)^{\ast }$, we obtain the commutativity condition for $Y(t)$ and $U(t,s)$%
, which is equivalent (\ref{eq:2.7}) for $X=Z_{\nu }^{\mu }$ due to the
approximation (\ref{eq:2.13}) of (\ref{eq:2.12}) and adaptedness of $Y$. The
condition (\ref{eq:2.7}) in the terms of $\widehat{Z}_{\nu }^{\mu }U=UZ_{\nu
}^{\mu }$ can be written as: 
\begin{equation*}
U(s)\left[ Y(t),Z_{\nu }^{\mu }(s)\right] =\left[ \widehat{Y}(t),\widehat{Z}%
_{\nu }^{\mu }(s)\right] U(s)=0,\forall t\leq s.
\end{equation*}%
In the same way the non-demolition condition for an output process $Y$ with
respect to a QS process $X$ can be written as $[\widehat{Y}(t),\widehat{X}%
(s)]U(s)=0$ in terms of $\widehat{X}U=UX$.

Representing $X$ in the case (\ref{eq:2.10}) in the form of (\ref{eq:2.5})
as the solution of the QS -- equation 
\begin{equation}
\mathrm{d}X=\left( \mathbf{Z}^{\star }\mathbf{FZ}-X\otimes \mathbf{\delta }%
\right) _{\nu }^{\mu }\mathrm{d}\widehat{A}_{\mu }^{\nu },\quad X(0)=%
\widehat{x}\otimes I,  \label{eq:2.14}
\end{equation}%
and taking into account that due to (\ref{eq:2.4}) 
\begin{equation}
\mathrm{d}(YX)=(\mathbf{Z}^{\star }\mathbf{GFZ}-YX\otimes \mathbf{\delta }%
)_{\nu }^{\mu }\mathrm{d}\widehat{A}_{\mu }^{\nu },  \label{eq:2.15}
\end{equation}%
one can easily obtain, that $[Y,X]=0$, iff $[\widehat{y},\widehat{x}]=0$ and 
$[\mathbf{G},\mathbf{F}]=0$ due to $\mathbf{ZZ}^{\star }\mathbf{F}=\mathbf{F}%
=\mathbf{FZZ}^{\star }$. In order to satisfy the condition $[Y(t),X(s)]=0$
for all $t\leq s$, it should be completed by $[Y(t),F_{\nu }^{\mu }(s)]=0$
at least for $\mu \not=+$ or $\nu \not=-$ and all $t\leq s$ due to the QS
integral representation 
\begin{equation*}
X(s)=X(t)+\int_{t}^{s}(\mathbf{Z}^{\star }\mathbf{FZ}-X\otimes \mathbf{%
\delta })_{\nu }^{\mu }\mathrm{d}\widehat{A}_{\mu }^{\nu }\;\;\mathrm{for}%
\;\;s>t
\end{equation*}%
and commutativity of $Y(t)$ with $\mathbf{Z}(s)$ at $s\geq t$. So $[\mathbf{D%
},\mathbf{F}]=[\mathbf{G},\mathbf{F}]-[Y\otimes \mathbf{\delta },\mathbf{F}%
]=0$, what gives the necessary and sufficient nondemolition conditions,
which can be written in the terms of $\widehat{Y},\widehat{\mathbf{D}},%
\widehat{\mathbf{F}}$ as (\ref{eq:2.11}) by multiplication on the right by
the corresponding $U$.\hfill \vrule height.9exwidth.8exdepth-.1ex

\begin{corollary}
The process $X$ is an evolute transformation $X=U^{\ast }(\widehat{x}\otimes
I)U$ of an initial operator $\widehat{x}\in {\mathcal{B}}(\mathfrak{h})$
with respect to a QS unitary process $U$, described by the QS equation (\ref%
{eq:2.6}) iff it satisfies the QS equation (\ref{eq:2.14}) with $\mathbf{F}%
=X\otimes \mathbf{\delta }$. The process $Y$ is an output process with
respect to the QS Markovian evolution defined on the von Neumann algebra ${%
\mathcal{A}}={\mathcal{B}}(\mathfrak{h})$ by the transformation $Z_{\nu
}^{\mu }=U^{\ast }\left( \widehat{z}_{\nu }^{\mu }\otimes I\right) U$ of the
initial QS generators $\widehat{z}_{\nu }^{\mu }$, acting in $\mathfrak{h}$,
iff $\left[ \widehat{Y}(t),\widehat{z}_{\nu }^{\mu }\otimes I\right] =0$ for
all $t$ and $\mu ,\nu $. The output process $\widehat{Y}$ is nondemolition
with respect to $X=U^{\ast }(\widehat{x}\otimes I)U$ for arbitrary $\widehat{%
x}\in {\mathcal{B}}(\mathfrak{h})$, iff $\widehat{Y}=\widehat{1}\otimes B$,
where $B$ is an adapted process in Fock space ${\mathcal{G}}$. \label%
{corol:2}
\end{corollary}

Indeed, if $\widehat{X}(t)=\widehat{x}\otimes I$ is a time independent
adapted process, then it satisfies the QS equation $\mathrm{d}\widehat{X}=%
\widehat{C}_{\nu }^{\mu }\mathrm{d}\widehat{A}_{\mu }^{\nu }$ corresponding
to $\widehat{C}_{\nu }^{\mu }=0=\widehat{F}_{\nu }^{\mu }-\widehat{z}\otimes
I\delta _{\nu }^{\mu }$. Hence, the process $\mathbf{Z}=U^{\ast }\widehat{%
\mathbf{Z}}U$ satisfies the equation (\ref{eq:2.14}) with $F_{\nu }^{\mu
}=U^{\ast }(\widehat{x}\otimes I)U\delta _{\nu }^{\mu }=X\delta _{\nu }^{\mu
}$.

The output condition $\left[ \widehat{Y}(t),\widehat{Z}_{\nu }^{\mu }(s)%
\right] U(s)=0$, for $\widehat{\mathbf{Z}}(s)=\widehat{\mathbf{z}}\otimes I$
and unitary $U$, means $\left[ \widehat{Y}(t),\widehat{z}_{\nu }^{\mu
}\otimes I\right] =0$ for all $t,\mu ,\nu $; moreover the nondemolition
condition $[\widehat{Y}(t),\widehat{X}(s)]U(s)=0$ for $\widehat{X}(s)=%
\widehat{z}\otimes I$ with arbitrary $\widehat{x}\in \mathcal{B}(\mathfrak{h}%
)$ is possible only if $Y=\widehat{1}\otimes B$.

Note that an output QS process $Y=U^{\ast }\widehat{Y}U$ is defined as the
sum of $\widehat{y}\otimes I$ and a QS -- integral 
\begin{equation*}
\int_{0}^{t}A(\mathbf{D},\mathrm{d}s)=U(t)^{\ast }\int_{0}^{t}\widehat{A}(%
\widehat{\mathbf{D}},\mathrm{d}s)U(t)
\end{equation*}%
with the QS differential $A(\mathbf{D},\mathrm{d}t)=(\mathbf{Z}^{\star }${$%
\mathbf{D}$}$\mathbf{Z})_{\nu }^{\mu }(t)\mathrm{d}\widehat{A}_{\mu }^{\nu
}(t)$. In the case of commuting matrix elements $D_{\mu }^{\kappa }=U^{\ast }%
\widehat{D}_{\mu }^{\kappa }U$ and $Z_{\nu }^{\iota }=U^{\ast }\widehat{Z}%
_{\nu }^{\iota }U$, as it happens for $\widehat{Z}_{\nu }^{\iota }(t)=%
\widehat{z}_{\nu }^{\iota }\otimes I$, $\widehat{D}_{\mu }^{\kappa }=%
\widehat{1}\otimes \widehat{D}_{\mu }^{\kappa }$, this integral can be
defined as the QS It\^{o} integral 
\begin{equation*}
\int_{0}^{t}A(\mathbf{D},\mathrm{d}s)=\int_{0}^{t}\left( D_{+}^{-}\mathrm{d}%
s+D_{\circ }^{-}\mathrm{d}A_{-}^{\circ }+D_{+}^{\circ }\mathrm{d}A_{\circ
}^{+}+D_{\circ }^{\circ }\mathrm{d}N_{\circ }^{\circ }\right)
=\int_{0}^{t}D_{\nu }^{\mu }\mathrm{d}A_{\mu }^{\nu }
\end{equation*}%
with respect to output annihilation $A_{-}^{\circ }$, creation $A_{\circ
}^{+}$ and quantum number $N_{\circ }^{\circ }$ processes 
\begin{equation*}
A_{-}^{\circ }(t)=\int_{0}^{t}\left( Z_{k}^{\circ }\mathrm{d}\widehat{A}%
_{-}^{k}+Z_{+}^{\circ }\mathrm{d}s\right) =A_{\circ }^{+}(t)^{\ast }\
,\qquad \qquad N_{\circ }^{\circ }(t)=
\end{equation*}%
\begin{equation*}
\int_{0}^{t}\left( (Z_{\circ }^{\circ ^{\ast }}Z_{\circ }^{\circ })_{k}^{i}%
\mathrm{d}\widehat{N}_{i}^{k}+(Z_{\circ }^{\circ ^{\ast }}Z_{+}^{\circ })^{i}%
\mathrm{d}\widehat{A}_{i}^{+}+(Z_{+}^{\circ ^{\ast }}Z_{\circ }^{\circ })_{k}%
\mathrm{d}\widehat{A}_{-}^{k}+(Z_{+}^{\circ \ast }Z_{+}^{\circ }\mathrm{d}%
s)\right)
\end{equation*}%
as the unitary transformation $A_{\nu }^{\mu }=U^{\ast }\widehat{A}_{\nu
}^{\mu }U$ of the input canonical processes $\widehat{A}_{\nu }^{\mu }$. 

\section{QS nonlinear nondemolition filtering}

\label{sec:belat3}

Let us consider a selfadjoint family $Y=(Y_{i})$ of commuting output
processes $Y_{i},i=1,\dots ,n$, defined by 
\begin{equation}
Y_{i}(t)=\widehat{y}_{i}\otimes I+\int_{0}^{t}(\mathbf{Z}^{\star }\mathbf{D}%
_{i}\mathbf{Z})_{\nu }^{\mu }\mathrm{d}\widehat{A}_{\mu }^{\nu },
\label{eq:3.1}
\end{equation}%
which are nondemolition with respect a QS process $X$: 
\begin{equation}
\left[ Y_{i}(t),Z_{\nu }^{\mu }(s)\right] =0,\;\left[ Y_{i}(t),X(s)\right]
=0\quad \forall t\leq s  \label{eq:3.2}
\end{equation}%
As follows from (\ref{eq:2.11}) for $\widehat{x}=\widehat{y}_{k}$, $\widehat{%
\mathbf{F}}=\widehat{\mathbf{D}}_{k}+\widehat{Y}_{k}\otimes \mathbf{\delta }$%
, the family $Y$ satisfies the selfnondemolition condition 
\begin{equation}
\left[ Y_{i}(t),Y_{k}(s)\right] =0,\quad \forall s,t;i,k,  \label{eq:3.3}
\end{equation}%
iff $\left[ \widehat{Y}_{i}(t),\widehat{D}_{\nu }^{\mu }(s)_{k}\right]
U(s)=0 $ for all $t\leq s$, and 
\begin{equation}
\left[ \widehat{y}_{i},\widehat{y}_{k}\right] =0,\;\;\left[ \widehat{\mathbf{%
D}}_{i},\widehat{\mathbf{D}}_{k}\right] U=0,\quad \forall i,k.
\label{eq:3.4}
\end{equation}

Let us denote by $\mathcal{A}_{t}=\left\{ Y_{i}^{t}:i=1,\dots ,n\right\}
^{\prime }$ the reduced algebra of bounded operators in $\mathcal{H}$,
corresponding to the measurements of the process $Y^{t}=\{Y(s):s\leq t\}$ up
to a time $t$, defined as the commutant of all $Y_{i}^{t}=U^{\ast }\widehat{Y%
}_{i}^{t}U$, and by $\mathcal{O}=\{y_{1},\dots y_{n}\}^{\prime }$ the
initial algebra, defining $\mathcal{A}_{0}=\mathcal{O}\otimes \mathcal{B}({%
\mathcal{F}})$. The nonincreasing family $(\mathcal{A}_{t})$ is the family
of maximal von Neumann subalgebras $\mathcal{A}_{s}\subseteq \mathcal{A}%
_{t}\subseteq \mathcal{A}_{0},\;s\geq t\geq 0$ of the initial reduced
algebra $\mathcal{A}_{0}$, with respect to which $Y$ is a nondemolition
commutative vector process in the sense of the definition $Y_{i}(t)\in 
\mathcal{A}_{t}^{\prime }\quad \forall t$ (or $Y_{i}(t)$ is affiliated to $%
\mathcal{A}_{t}^{\prime }$) of a nondemolition QS process given in \cite%
{bib:3}. The Abelian algebra $\mathcal{B}^{t}=\mathcal{A}_{t}^{\prime }$
generated by $Y^{t}$, with $\mathcal{B}^{0}=\mathcal{O}^{\prime }\otimes I$,
generated by $Y^{0}=y\otimes I$, forms the center $\mathcal{B}^{t}=\mathcal{A%
}_{t}\bigcap \mathcal{A}_{t}^{\prime }$ of $\mathcal{A}_{t}$, hence $%
\mathcal{A}_{t}$ is a decomposable algebra, having the conditional
expectations with respect to $\mathcal{B}^{t}=\mathcal{A}_{t}$ for any
normal initial state on $\mathcal{A}_{0}\supseteq \mathcal{A}_{t}$.

As follows from the next theorem, the nondemolition principle is not only
sufficient, but also necessary for the existence of compatible conditional
expectation on $\mathcal{A}_t$ with respect to $\mathcal{B}_t$ for an
arbitrary initial state vector $\xi$.

We shall explicitly construct the conditional expectation not only for
bounded $X\in \mathcal{A}_{t}$, but also for $X$ affiliated to $\mathcal{A}%
_{t}$. An operator $B$ is said to be defined almost everywhere with respect
to the pair $(\mathcal{B}^{t},\xi )$, if it is densely defined in the
support subspace $\mathcal{K}^{t}=P^{t}\mathcal{H}$, where $P^{t}=\inf
\{P=P^{\ast }P\in \mathcal{B}^{t}:P\xi =\xi \}$, $\Vert \xi \Vert =1$.

\begin{theorem}
Let $\mathcal{B}^{t}\subset \mathcal{A}_{t}$ be a von Neumann subalgebra on
a Hilbert space $\mathcal{H}$. Then a conditional expectation $\epsilon _{t}$
as a positive projection onto $\mathcal{B}^{t}$, satisfying the
compatibility condition $<\xi |\epsilon _{t}(X)\xi >=<\xi |X\xi >$, for all $%
X\in \mathcal{A}$ exists on $\mathcal{A}_{t}$ for an arbitrary $\xi \in 
\mathcal{H}$, iff $\mathcal{B}^{t}$ commutes with $\mathcal{A}_{t}:\mathcal{B%
}^{t}\subseteq \mathcal{A}_{t}^{\prime }$. In this case the algebra $%
\mathcal{A}_{t}$ can be extended to the commutant of $\mathcal{B}^{t}$, such
that for any operator $X$, commuting on the domain $\mathcal{A}_{t}\xi $
with $\mathcal{B}^{t}=\mathcal{A}_{t}^{\prime }$, the expectation $\epsilon
_{t}(X)$ is given on $\mathcal{A}_{t}\xi $ by 
\begin{equation}
\epsilon _{t}(X)A\xi =AE_{t}X\xi ,\quad \forall A\in \mathcal{A}_{t},
\label{eq:3.5}
\end{equation}%
where $E_{t}\in \mathcal{A}_{t}$ is the orthoprojector on ${\overline{%
\mathcal{B}^{t}\xi }}$. The formula (\ref{eq:3.5}) uniquely defines $%
\epsilon _{t}(X)$ as an operator $\epsilon _{t}(X)P^{t}$ affiliated to $%
\mathcal{B}^{t}$ on $\overline{\mathcal{A}_{t}\xi }$ even for unbounded $X$. %
\label{th:3}
\end{theorem}

\noindent \textsc{Proof.}\ Let us suppose that $[X,B]\not=0$ for an $X\in 
\mathcal{A}_{t}$ and $B\in \mathcal{B}^{t}$, and that $\epsilon _{t}:%
\mathcal{A}_{t}\rightarrow \mathcal{B}^{t}$ is defined as a positive
projection, compatible with $\xi \in \mathcal{H}$, for which $<\xi |[X,B]\xi
>\not=0$. Then, due to the modularity property 
\begin{equation*}
\epsilon _{t}(XB)=\epsilon _{t}(X)B\ ,\;\epsilon _{t}(BX)=B\epsilon _{t}(X)\
,
\end{equation*}%
where $X\in \mathcal{A}_{t}$, $B\in \mathcal{B}^{t}$, we would have $<\xi
|[\epsilon _{t}(x),B]\xi >=<\xi |[X,B]\xi >\not=0$, what would be possible
only if $\mathcal{B}^{t}$ would be non-Abelian. But for non-Abelian $%
\mathcal{B}^{t}$ the conditional expectation does not exist for all vectors $%
\xi \in \mathcal{H}$, as can be easily shown for a factor $\mathcal{B}%
^{t}\not=\mathbb{C}\widehat{I}$. Indeed, in this case such a vector $\xi $
has to be of the form $\xi _{0}\otimes \xi _{1}$, and $\epsilon (A\otimes
B)=<\xi _{0}|A\xi _{0}>I_{0}\otimes B$, where $\xi _{0}\in \mathcal{H}%
_{0}\not=\mathbb{C}$, if $\mathcal{A}_{t}\not=\mathcal{B}^{t}$, $\xi _{1}\in 
\mathcal{H}_{1}=\overline{\mathcal{B}^{t}\xi }$, corresponding to the
decomposition $\mathcal{H}=\mathcal{H}_{0}\otimes \mathcal{H}_{1}$. So, it
is necessary that $\mathcal{B}^{t}\subseteq \mathcal{A}_{t}^{\prime }$.

Let us define $\epsilon _{t}$ for such an Abelian algebra $\mathcal{B}^{t}$
by (\ref{eq:3.5}) with $\mathcal{A}_{t}^{\prime }=\mathcal{B}^{t}$ and a
fixed $\xi \in \mathcal{H}$. The orthoprojector $E_{t}$ commutes with $%
\mathcal{A}_{t}^{\prime }$ due to the invariance of $\mathcal{E}_{t}=%
\overline{\mathcal{A}_{t}^{\prime }\xi }$ with respect to the action of the
algebra $\mathcal{A}_{t}^{\prime }$. Hence the operator $E_{t}XE_{t}$
commutes with $\mathcal{A}_{t}^{\prime }E_{t}$: 
\begin{equation*}
E_{t}XE_{t}BE_{t}=E_{t}XBE_{t}=E_{t}BXE_{t}=BE_{t}XE_{t}=E_{t}BE_{t}XE_{t}\ ,
\end{equation*}%
if the operator $X$ commutes with the all $B\in \mathcal{A}_{t}^{\prime }$.
But this means that $E_{t}XE_{t}$ is affiliated with the reduced von Neumann
algebra $E_{t}\mathcal{A}_{t}E_{t}$ on $\mathcal{E}_{t}$, coinciding with
its commutant $\mathcal{A}_{t}^{\prime }E_{t}$ on $\mathcal{E}_{t}$ because
the induced Abelian algebra $\mathcal{A}_{t}^{\prime }E_{t}$ has the cyclic
vector $\xi $ in $\mathcal{E}_{t}$. The commutativity of $E_{t}\mathcal{A}%
_{t}E_{t}=\mathcal{A}_{t}^{\prime }E_{t}$ helps to establish the correctness
of the definition (\ref{eq:3.5}) of the linear operator $\epsilon _{t}(X)$
on $\mathcal{A}_{t}\xi $ $A\xi =0\Rightarrow \epsilon _{t}(X)A\xi =0$.
Indeed, $\Vert \epsilon _{t}(x)A\xi \Vert =$ 
\begin{equation*}
\Vert AE_{t}X\xi \Vert =\Vert (E_{t}A^{\ast }AE_{t})^{1/2}E_{t}XE_{t}\xi
\Vert =\Vert E_{t}XE_{t}(E_{t}A^{\ast }AE_{t})^{1/2}\xi \Vert
\end{equation*}%
because $E_{t}\xi =\xi $ and $(E_{t}A^{\ast }AE_{t})^{1/2}$ $\xi =0$, if $%
A\xi =0$.

The operator $\epsilon _{t}(X):A\xi \rightarrow AE_{t}X\xi $ having the
range $\mathcal{A}_{t}E_{t}X\xi \subseteq \mathcal{K}^{t}=\overline{\mathcal{%
A}_{t}\xi }$, commutes with arbitrary $A\in \mathcal{A}_{t}$ due to the
definition (\ref{eq:3.5}), so $P^{t}\epsilon _{t}(X)$ is affiliated to $P^{t}%
\mathcal{A}_{t}^{\prime }P^{t}$, coinciding with $\mathcal{A}_{t}^{\prime
}P^{t}$ because $P^{t}\in \mathcal{A}_{t}^{\prime }\bigcap \mathcal{A}_{t}=%
\mathcal{A}_{t}^{\prime }$, if $P^{t}$ is the orthoprojector on $\mathcal{K}%
^{t}$.

The map $X\mapsto \epsilon _{t}(X)$ satisfies the unital property $\epsilon
^{t}(\widehat{I})A\xi =AE_{t}\xi =A\xi $ due to $\xi \in \mathcal{E}_{t}$,
and the modularity property 
\begin{equation*}
\epsilon _{t}(XB)A\xi =AE_{t}BX\xi =ABE_{t}X\xi =BAE_{t}X\xi =B\epsilon
_{t}(X)A\xi
\end{equation*}%
for all $A\in \mathcal{A}_{t}$ and $B\in \mathcal{A}_{t}^{\prime }$, and,
hence, maps the algebra $\mathcal{A}_{t}$ on the subalgebra $\mathcal{A}%
_{t}^{\prime }\subset \mathcal{A}_{t}$, represented on $\mathcal{K}^{t}$.

Now let us prove the uniqueness of the representation (\ref{eq:3.5}) of
conditional expectation $\epsilon _{t}$ as a map onto factor subalgebra $%
\mathcal{A}_{t}^{\prime }/\mathcal{A}_{t}^{\prime }P_{1}^{t}=\mathcal{A}%
_{t}^{\prime }P^{t}$, where $P^{t}=\widehat{I}-P_{1}^{t}\in \mathcal{A}%
_{t}^{\prime }$ is the support of $\xi $ which is the orthoprojector on $%
\overline{\mathcal{A}_{t}\xi }=\mathcal{K}^{t}$. Due to the commutativity of 
$\epsilon (X)$ with $\mathcal{A}_{t}$ we have $\epsilon _{t}(X)A\xi
=A\epsilon _{t}(X)\xi $ for $A\in \mathcal{A}_{t}$. So we have to prove,
that $\epsilon _{t}(X)\xi =E_{t}X\xi $. But $\epsilon _{t}(X)\xi \in 
\mathcal{A}_{t}^{\prime }\xi $, because $\epsilon _{t}(X)\in \mathcal{A}%
_{t}^{\prime }$ for $X\in \mathcal{A}_{t}$; hence we should prove, that $%
<B\xi |\epsilon _{t}(X)\xi >=<B\xi |E_{t}X\xi >$ for all $B\in \mathcal{A}%
_{t}^{\prime }$, which is a consequence of modularity and compatibility
conditions: 
\begin{equation*}
<B\xi |\epsilon _{t}(X)\xi >=<\xi |\epsilon _{t}(B^{\ast }X)\xi >=<\xi
|B^{\ast }X\xi >=<B\xi |X\xi >=<B\xi |E_{t}X\xi >\ .
\end{equation*}%
\medskip

\noindent \textbf{Remark.\/} Note that one should identify the
factor-algebra $\mathcal{B}^{t}P^{t}$ with the space $L^{\infty }(\mathcal{V}%
^{t})$ of essentially bounded measurable complex functions on the
probability space $\mathcal{V}^{t}$ of all observed values $%
v^{t}=\{v(s):s\leq t\}$, $v(t)=(v_{i})(t)$ of the commutative vector process 
$Y^{t}$, stopped at $t$. The probability measure $\mu (\mathrm{d}v^{t})=<\xi
|I(\mathrm{d}v^{t})\xi >$ is induced on the Borel $\sigma $-algebra of $%
\mathcal{V}^{t}$ by the spectral resolution $Y^{t}=\int v^{t}I(\mathrm{d}%
v^{t})$. If $P^{t}=\int_{\mathcal{V}^{t}}^{\otimes }P_{v^{t}}\mu (\mathrm{d}%
v^{t})$ is the corresponding decomposition of $P^{t}\in \mathcal{A}%
_{t}^{\prime }$, then 
\begin{equation}
P^{t}\epsilon _{t}(X)=\int_{\mathcal{V}^{t}}^{\otimes
}<X>_{v^{t}}P_{v^{t}}\mu (\mathrm{d}v^{t})\ ,  \label{eq:3.6}
\end{equation}%
where $<X>_{v^{t}}=<\xi _{v^{t}}|X\xi _{v^{t}}>$, and the vectors $\xi
_{v^{t}}=P_{v^{t}}\xi /\Vert P_{v^{t}}\xi \Vert $ define the resolution $%
E_{t}=\int_{\mathcal{V}^{t}}^{\otimes }|\xi _{v^{t}}><\xi _{v^{t}}|\mu (%
\mathrm{d}v^{t})$. Hence one should consider $\epsilon _{t}(X)$ for an $X\in 
\mathcal{A}_{t}$ as a function $\epsilon _{t}(X):\mathcal{V}^{t}\rightarrow 
\mathbb{C}P_{v^{t}}$ giving for almost all trajectories $v^{t}\in \mathcal{V}%
^{t}$, observed up to a time $t$, the posterior mean values $<X>_{v^{t}}$ of
a QS nondemolished process $X(t)$. The initial conditional expectation $%
\epsilon _{0}$ with respect to $\mathcal{B}^{0}=\mathcal{O}^{\prime }\otimes
I$ and $\xi =\psi \otimes \varphi $ is given for $X=\widehat{x}\otimes I$ as 
$\epsilon (\widehat{X})\otimes I$ by 
\begin{equation}
\widehat{p}\epsilon (\widehat{x})=\int_{\mathcal{V}_{0}}^{\oplus }<\psi _{v}|%
\widehat{x}\psi _{v}>\widehat{p}_{v}\mu (\mathrm{d}v)\ ,\quad \mu (\mathrm{d}%
v)=\Vert \widehat{1}(\mathrm{d}v)\psi \Vert ^{2}\ .  \label{eq:3.7}
\end{equation}%
Here the vectors $\psi _{v}=\widehat{p}_{v}\psi /\Vert \widehat{\psi }%
_{v}\psi \Vert $, $v\in \mathcal{V}_{0}$ and the decomposition $\widehat{e}%
=\int_{\mathcal{V}}^{\otimes }|\psi _{v}><\psi _{v}|\mu (\mathrm{d}v)$ for $%
E_{0}=\widehat{e}\otimes I$, $\{p_{v}\}$ define the decomposition $\widehat{p%
}=\int^{\otimes }\widehat{p}_{v}\mu (\mathrm{d}v)$ for $P^{0}=\widehat{p}%
\otimes I$, corresponding to the orthogonal resolution $\widehat{y}=\int v%
\widehat{1}(\mathrm{d}v)$ on the spectrum $\mathcal{V}_{0}$ of the
commutative family $\widehat{y}=(\widehat{y}_{i})$ of the initial operators $%
Y(0)=\widehat{y}\otimes I$. 

\section{QS calculus of a posteriori expectations}

\label{sec:belat4}

Now let us suppose that $\xi =\psi \otimes \delta _{\emptyset }$, the output
commuting processes (\ref{eq:3.1}) are nondemolition with respect to 
\begin{equation}
X(t)=\widehat{x}\otimes I+\int_{0}^{t}(\mathbf{Z}^{\star }\mathbf{FZ}%
-X\otimes \mathbf{\delta })_{\nu }^{\mu }\mathrm{d}\widehat{A}_{\mu }^{\nu }
\label{eq:4.1}
\end{equation}%
with $\widehat{x}\in \mathcal{O}$, and $\mathbf{F}(t)\in \mathcal{F}_{t}$
where $\mathcal{F}_{t}$ is the $\star $-algebra of matrix-operators $\mathbf{%
F}=(F_{\nu }^{\mu })$, commuting with $Y(s)$, $s\leq t$ and $\mathbf{D}(t):%
\mathcal{F}_{t}=\{F_{\nu }^{\mu }\in \mathcal{A}_{t}:[\mathbf{D}_{i},\mathbf{%
F}]=0$, $i=1,\dots ,n\}$. In the following we shall also demand that the
process $Y=(Y_{i})$ is continuous from the right in the sense $\mathcal{A}%
_{t}^{\prime }=\bigcap_{s>t}\mathcal{A}_{s}^{\prime }$, what is equivalent
to $D_{i}(t)_{\nu }^{\mu }\in \mathcal{A}_{t}^{\prime }$ for all $i$ and $t$.

Let us denote by $\mathcal{C}$ the linear span of the initial operators $%
\{y_{i}\}$ with the operators $y_{0}\in \mathcal{C}_{0}$ from the
commutative ideal $\mathcal{C}_{0}=\{b\in \mathcal{O}^{\prime }:<b\psi
|b\psi >=0\}$. We also denote by $\mathcal{D}^{t}$ the $\mathcal{A}%
_{t}^{\prime }$-span of the operator-matrices $\{\mathbf{D}_{i}\}(t)$ with
the ideal%
\begin{equation*}
\mathcal{D}_{0}^{t}=\{\mathbf{D}\in \mathcal{F}_{t}^{\prime
}:D_{-}^{-}=0=D_{+}^{+},(\mathbf{D}Z_{+}|\mathbf{D}Z_{+})=0\}
\end{equation*}%
of the commutative $\star $-algebra%
\begin{equation*}
\mathcal{F}_{t}^{\prime }=\{D_{\nu }^{\mu }\in \mathcal{A}_{t}^{\prime }|[%
\mathbf{D},\mathbf{F}]=0,\mathbf{F}\in \mathcal{F}\},
\end{equation*}%
corresponding to the kernel of the pseudoscalar product 
\begin{equation*}
(\mathbf{Z}_{+}|\mathbf{Z}_{+})=<\xi |(\mathbf{Z}^{\star }\mathbf{Z}%
)_{+}^{-}\xi >=<\xi |\mathbf{Z}_{+}^{\ast }\mathbf{gZ}_{+}\xi >\ ,
\end{equation*}%
where $\mathbf{Z}_{+}^{\ast }=\left( Z_{+}^{-\ast },Z_{+}^{\circ \ast },%
\widehat{I}\right) $ is the conjugate row to the column-operator $\mathbf{Z}%
_{+}$. Now we can formulate the main theorem. 

\begin{theorem}
Suppose that the output observed process (\ref{eq:3.1}) is nondemolition
with respect to a QS process (\ref{eq:4.1}), and the spans $\mathcal{C}%
\subseteq \mathcal{O}^{\prime }$ and $\mathcal{D}^{t}\subseteq \mathcal{F}%
_{t}^{\prime }$ are $^{\ast }$ -- and $\star $ -- algebras correspondingly.
Then the posterior mean value $\epsilon _{t}(X(t))$ for an initial state
vector $\xi =\psi \otimes \delta _{\emptyset },\psi \in \mathfrak{h}$, is
defined by an adapted commutative vector--process $\kappa _{t}=\left( \kappa
_{t}^{i}\right) ,\kappa _{t}^{i}\in \mathcal{A}_{t}^{\prime },i=1,\dots ,n$,
almost everywhere as an $\mathcal{A}_{t}^{\prime }$ linear nonaticipating
transformation of the output process $Y$ by the stochastic It\^{o} equation 
\begin{equation}
\mathrm{d}\epsilon _{t}(X(t))=\epsilon _{t}(\mathbf{Z}^{\star }\mathbf{F}%
\mathbf{Z})_{+}^{-}(t)\mathrm{d}t+\kappa _{t}^{i}(X(t))\mathrm{d}\widetilde{Y%
}_{i}(t)  \label{eq:4.2}
\end{equation}%
Here $\epsilon _{t}(\mathbf{F})_{+}^{-}=\epsilon _{t}(F_{+}^{-}),\epsilon
_{0}(\widehat{x}\otimes I)=\epsilon (\widehat{x})\otimes I$, $\kappa ^{i}%
\mathrm{d}Y_{i}\equiv \sum_{i=1}^{n}\kappa ^{i}\mathrm{d}Y_{i}$, and 
\begin{equation}
\mathrm{d}\widetilde{Y}_{i}(t)=\mathrm{d}Y_{i}(t)-\epsilon _{t}(\mathbf{Z}%
^{\star }\mathbf{D}_{i}\mathbf{Z})_{+}^{-}(t)\mathrm{d}t,\;\widetilde{Y}%
_{i}(0)=\widetilde{y}_{i}\otimes I  \label{eq:4.3}
\end{equation}%
are the observed martingales with respect to the filtration $\left( {%
\epsilon }_{t}\right) $, and state vector $\xi $, called the innovating
process for $\left( \mathcal{A}_{t}^{\prime }\right) $. The process $\kappa
_{t}$ is defined uniquely up to the kernel of the correlation
matrix--process 
\begin{equation}
\varrho _{ik}(t)=\epsilon _{t}\left( \mathbf{Z}^{\star }\mathbf{D}%
_{i}^{\star }\mathbf{D}_{k}\mathbf{Z}\right) _{+}^{-}(t)=\epsilon
_{t}[(D_{\circ }^{\circ }Z_{\circ }^{\circ }+D_{+}^{\circ })_{i}^{\ast
}(D_{\circ }^{\circ }Z_{+}^{\circ }+D_{+}^{\circ })_{k}](t)  \label{eq:4.4}
\end{equation}%
by the linear algebraic equation 
\begin{equation}
\varrho _{ik}(t)\kappa _{t}^{k}=\epsilon _{t}\left( \mathbf{Z}^{\star }%
\mathbf{D}_{i}^{\star }\mathbf{F}\mathbf{Z}\right) _{+}^{-}(t)-\epsilon
_{t}(X(t))\epsilon _{t}\left( \mathbf{Z}^{\star }\mathbf{D}_{i}^{\star }%
\mathbf{Z}\right) _{+}^{-}(t),  \label{eq:4.5}
\end{equation}%
having in the case $\mathbf{F}=X\otimes \mathbf{\delta }$, corresponding to $%
\widehat{X}(t)=\widehat{x}\otimes I$, the form 
\begin{equation}
\varrho _{ik}(t)\kappa _{t}^{k}=\epsilon _{t}\left( Z_{+}^{\ast }D_{\circ
i}^{\circ \ast }\tilde{X}Z_{+}^{\circ }\right) (t)+\epsilon _{t}\left( 
\widetilde{X}D_{+i}^{\circ \ast }Z_{+}^{\circ }+Z_{+}^{\circ \ast }D_{\circ
i}^{-\ast }\widetilde{X}\right) (t),  \label{eq:4.6}
\end{equation}%
where $\widetilde{X}(t)=X(t)-\epsilon _{t}(X(t))$. The initial a posteriori
mean value $\epsilon (\widehat{x})$ is the linear combination $\epsilon (%
\widehat{x})=<\psi |\widehat{x}\psi >+\kappa ^{i}(\widehat{x})\widetilde{y}%
_{i}$, of $\widetilde{y}_{i}=\widehat{y}_{i}-<\psi |\widehat{y}_{i}\psi >%
\widehat{1}$, where $\kappa =(\kappa ^{i})$ is defined by the equation 
\begin{equation}
\varrho _{ik}\kappa ^{k}=<\psi |\widehat{y}_{i}^{\ast }\widehat{x}\psi
>,\quad \varrho _{ik}=<\psi |\widetilde{y}_{i}^{\ast }\widetilde{y}_{k}\psi
>,  \label{eq:4.7}
\end{equation}%
with $\widetilde{x}=\widetilde{x}-<\psi |\widehat{x}\psi >\widehat{1}$,
uniquely up to the kernel of the initial correlation matrix $\varrho
=(\varrho _{ik})$. \label{th:4}
\end{theorem}

In order to prove this fundamental filtering theorem we need the following
lemmas.

\begin{lemma}
If the process $X$ satisfies the equation (\ref{eq:4.1}), then there exists
such a martingale $M_{t}$ with respect to $\left( \epsilon _{t},\xi \right) $%
, affiliated with $\mathcal{A}_{t}^{\prime }$ on $\mathcal{A}_{t}\xi $, such
that almost everywhere 
\begin{equation}
\epsilon _{t}(X(t))=\epsilon (\widehat{x})\otimes I+\int_{0}^{t}\epsilon
_{s}(\mathbf{Z}^{\star }\mathbf{F}\mathbf{Z})_{+}^{-}(s)\mathrm{d}s+M_{t}.
\label{eq:4.8}
\end{equation}%
\label{lemma:4.1}
\end{lemma}

\noindent \textsc{Proof.}\ Let us define $M_{t}$ on $\xi $ by 
\begin{equation*}
M_{t}\xi =\left( E_{t}-E_{0}\right) X(0)\xi +\int_{0}^{t}\left(
E_{t}-E_{s}\right) (\mathbf{Z}^{\star }\mathbf{FZ})_{\nu }^{\mu }\mathrm{d}%
A_{\mu }^{\nu }\xi \ .
\end{equation*}%
Obviously, that $M_{t}\xi $ satisfies the $\left( \epsilon _{t},\xi \right) $
-- martingale condition $E_{s}M_{t}\xi =M_{s}\xi $ for all $s\leq t$, and 
\begin{equation*}
E_{t}X(t)\xi =\left( \widehat{e}\widehat{x}\otimes I\right) \xi
+\int_{0}^{t}E_{s}(\mathbf{Z}^{\star }\mathbf{FZ})_{+}^{-}(s)\xi \mathrm{d}%
s+M_{t}\xi
\end{equation*}%
due to $E_{s}(\mathbf{Z}^{\star }\mathbf{FZ})_{\nu }^{\mu }\mathrm{d}A_{\mu
}^{\nu }\xi =E_{s}(\mathbf{Z}^{\star }\mathbf{FZ})_{+}^{-}\xi \mathrm{d}t$
for $\xi =\psi \otimes \delta _{\emptyset }$. The operator $M_{t}$,
affiliated with $\mathcal{A}_{t}^{\prime }$ can be correctly defined almost
everywhere by 
\begin{equation*}
M_{t}A\xi =AE_{t}M_{t}\xi =AM_{t}\xi ,\quad \forall A\in \mathcal{A}_{t},
\end{equation*}%
as in the case of (\ref{eq:3.5}) for $X=M_{t},\;\epsilon _{t}(M_{t})=M_{t}$.

So, for any $A\in \mathcal{A}_{t}$ we have 
\begin{equation*}
\epsilon _{t}(X(t))A\xi =A\left( \widehat{e}\widehat{x}\psi \otimes \delta
_{\emptyset }+\int_{0}^{t}E_{s}(\mathbf{Z}^{\star }\mathbf{FZ}%
)_{+}^{-}(s)\xi \mathrm{d}s+M_{t}\xi \right) =
\end{equation*}%
\begin{equation*}
=\left( \epsilon (\widehat{x})\otimes I\right) A\xi +\int_{0}^{t}\epsilon
_{s}(\mathbf{Z}^{\star }\mathbf{FZ})_{+}^{-}(s)A\xi \mathrm{d}s+M_{t}A\xi ,
\end{equation*}%
and, hence, (\ref{eq:4.8}) holds on the dense linear manifold $\mathcal{A}%
_{t}\xi $ of the support $\mathcal{K}^{t}$ of the state $\xi $ on $\mathcal{A%
}_{t}^{\prime }$. \medskip

\begin{lemma}
A process $M_{t}=\int_{0}^{t}(\mathbf{Z}^{\star }\mathbf{D}\mathbf{Z})_{\nu
}^{\mu }\mathrm{d}\widehat{A}_{\mu }^{\nu }$ with $\mathbf{D}(t)\in \mathcal{%
D}^{t}$ is a martingale with respect to $(\epsilon _{t},\xi )$, iff $%
\epsilon _{t}(\mathbf{Z}^{\star }\mathbf{D}\mathbf{Z})_{+}^{-}(t)=0$ for all 
$t$, that is almost everywhere 
\begin{equation}
D_{+}^{-}(t)+D_{\circ }^{-}(t)\epsilon _{t}(Z_{+}^{\circ })+\epsilon
_{t}(Z_{+}^{\circ })^{\ast }D_{+}^{\circ }(t)+\epsilon _{t}(Z_{+}^{\circ
\ast }D_{\circ }^{\circ }Z_{+}^{\circ })(t)=0  \label{eq:4.9}
\end{equation}%
and is the zero martingale (almost everywhere), iff $\mathbf{D}(t)\in 
\mathcal{D}_{0}^{t}$, which is equivalent to $\epsilon _{t}(\mathbf{Z}%
^{\star }\mathbf{D}^{\star }\mathbf{D}\mathbf{Z})_{+}^{-}(t)=0$ for all $t$
almost everywhere, that is 
\begin{equation}
\epsilon _{t}[(D_{+}^{\circ }+D_{\circ }^{\circ }Z_{+}^{\circ })^{\ast
}(D_{+}^{\circ }+D_{\circ }^{\circ }Z_{+}^{\circ })](t)=0  \label{eq:4.10}
\end{equation}%
and, hence, $D_{+}^{-}(t)=\epsilon _{t}(Z_{+}^{\circ \ast }D_{\circ }^{\circ
}Z_{+}^{\circ })(t)$. 
\label{lemma:4.2}
\end{lemma}

\noindent \textsc{Proof.\/} Due to commutativity of $M(t)$ with $\mathcal{A}%
_{t}$, we have to prove only that $E_{t}M_{r}\xi =M_{t}\xi $ for all $r\geq
t $, iff $E_{t}(\mathbf{Z}^{\star }\mathbf{DZ})_{+}^{-}(t)=0$ for all $t$.
Indeed, $E_{t}(M_{r}-M_{t})\xi =\int_{t}^{r}E_{t}(\mathbf{Z}^{\star }\mathbf{%
DZ})_{+}^{-}(s)\xi \mathrm{d}s=0$ for all $r>t$ iff $E_{t}(\mathbf{Z}^{\star
}\mathbf{DZ})_{+}^{-}(s)\xi =0$ for all $t\leq s$, which is equivalent to $%
E_{t}(\mathbf{Z}^{\star }\mathbf{DZ})_{+}^{-}(t)\xi =0$ for all $t$ due to $%
E_{t}E_{s}=E_{t}$ for $t\leq s$, written in the form (\ref{eq:4.10}) for 
\begin{equation*}
(\mathbf{Z}^{\star }\mathbf{DZ})_{+}^{-}=D_{+}^{-}+D_{\circ
}^{-}Z_{+}^{\circ }+Z_{+}^{\circ \ast }D_{+}^{\circ }+Z_{+}^{\circ \ast
}D_{\circ }^{\circ }Z_{+}^{\circ }\ .
\end{equation*}%
If $M_{t}$ is a martingale, then 
\begin{equation*}
\epsilon _{t}[(M_{r}-M_{t})^{\ast }(M_{r}-M_{t})]=\epsilon _{t}(M_{r}^{\ast
}M_{r})-M_{t}^{\ast }M_{t}=
\end{equation*}%
\begin{equation*}
\int_{t}^{r}\epsilon _{t}(\mathbf{Z}^{\star }\mathbf{D}^{\star }\mathbf{DZ}%
)_{+}^{-}(s)\mathrm{d}s\geq 0\ .
\end{equation*}%
Hence, if $M_{t}$ is a zero martingale, $|M_{r}-M_{t}|^{2}=0$, and $\epsilon
_{t}(\mathbf{Z}^{\star }\mathbf{D}^{\star }\mathbf{DZ})_{+}^{-}(s)=0$ for
all $t\geq s$, what is equivalent to $\epsilon _{t}(\mathbf{Z}^{\star }%
\mathbf{D}^{\star }\mathbf{DZ})_{+}^{-}(t)=0$ for all $t$, or to (\ref%
{eq:4.10}) in view of 
\begin{equation*}
(\mathbf{Z}^{\star }\mathbf{D}^{\star }\mathbf{DZ})_{+}^{-}=(D_{+}^{\circ
}+D_{\circ }^{\circ }Z_{+}^{\circ })^{\ast }(D_{+}^{\circ }+D_{\circ
}^{\circ }Z_{+}^{\circ })\ .
\end{equation*}%
But this means, that $(D_{+}^{\circ }+D_{\circ }^{\circ }Z_{+}^{\circ })\xi
=0$, i.e. $\mathbf{D}(t)\in \mathcal{D}_{0}^{t}$. Conversely if $\mathbf{D}%
(t)\in \mathcal{D}_{0}^{t}$, i.e. if $\langle \xi |(\mathbf{Z}^{\star }%
\mathbf{D}^{\star }\mathbf{DZ})_{+}^{-}\xi \rangle =0$, then $E_{t}(\mathbf{Z%
}^{\star }\mathbf{D}^{\star }\mathbf{DZ})_{+}^{-}(t)\xi =0$ because 
\begin{equation*}
\langle B\xi |(\mathbf{Z}^{\star }\mathbf{D}^{\star }\mathbf{DZ})_{+}^{-}\xi
\rangle =\langle \xi |(\mathbf{Z}^{\star }\mathbf{D}^{\star }(B\otimes 
\mathbf{\delta })\mathbf{DZ})_{+}^{-}\xi \rangle =0
\end{equation*}%
for any $B\in \mathcal{A}_{t}^{\prime }$. 

\begin{lemma}
Let the linear complex span of $\{\widehat{y}_{i}\}$ and the span of $\{%
\mathbf{D}_{i}\}(t)$ with the coefficients in $\mathcal{A}_{t}^{\prime }$,
be commutative $^{\ast }$ -- and $\star $ -- algebras $\mathcal{C}$ and $%
\mathcal{D}^{t}$ up to the ideals $\mathcal{C}_{\circ }\subseteq \mathcal{C}$
and $\mathcal{D}_{0}^{t}\subset \mathcal{D}^{t}$ correspondingly. Then the
locally bounded process 
\begin{equation}
B(t)=(\widehat{y}_{0}+\lambda _{0}^{i}\tilde{y}_{i})\otimes I+\int_{0}^{t}(%
\mathbf{Z}^{\star }(\mathbf{D}_{0}+\lambda _{s}^{i}\tilde{\mathbf{D}}_{i})%
\mathbf{Z})_{\nu }^{\mu }(s)\mathrm{d}A_{\mu }^{\nu }(s)\ ,  \label{eq:4.11}
\end{equation}%
where $y_{0}\in \mathcal{C}$, $\tilde{y}_{i}=\widehat{y}_{i}-<\psi |\widehat{%
y}_{i}\psi >\widehat{1}$, $D_{0}\in \mathcal{D}_{0}^{t}$, $\tilde{D}_{\nu
}^{\mu }=D_{\nu }^{\mu }$, if $(\mu ,\nu )\not=(-,+)$, and $\tilde{D}%
_{+}^{-}(t)=D_{+}^{-}(t)-\epsilon _{t}(\mathbf{Z}^{\star }\mathbf{D}\mathbf{Z%
})_{+}^{-}(t)$, defined by weakly measurable locally bounded functions $%
t\mapsto \lambda _{t}^{i}\in \mathcal{A}_{t}^{\prime }$, $\lambda
_{0}^{i}\in \mathbb{C}$, compose a weakly dense $^{\ast }$-algebra $\mathcal{%
C}^{t}$ in $\mathcal{A}_{t}^{\prime }$. 
\label{lemma:4.3}
\end{lemma}

\noindent \textsc{Proof.\/} Using the QS It\^{o} formula (\ref{eq:2.3}), one
obtains for $\mathrm{d}B=(\mathbf{Z}^{\star }\mathbf{DZ})_{\nu }^{\mu }%
\mathrm{d}A_{\mu }^{\nu }$ with $\mathbf{D}(t)=\mathbf{D}_{0}(t)+\lambda
_{t}^{i}\mathbf{\hat{D}}_{i}(t)\in \mathcal{D}^{t}$ 
\begin{equation*}
\mathrm{d}(B^{\ast }B)=(\mathbf{Z}^{\star }(\mathbf{D}^{\star }(B\otimes 
\mathbf{\delta })+(B\otimes \mathbf{\delta })^{\ast }\mathbf{D}+\mathbf{D}%
^{\star }\mathbf{D})\mathbf{Z})_{\nu }^{\mu }\mathrm{d}\widehat{A}_{\mu
}^{\nu }
\end{equation*}%
due to the commutativity of $B(t)$ with $Z_{\nu }^{\mu }(t)$. But $\mathbf{D}%
^{\star }(t)\mathbf{D}(t)\in \mathcal{D}^{t}$ and, hence $(\mathbf{D}^{\star
}(B\otimes \mathbf{\delta })+(B\otimes \mathbf{\delta })^{\ast }\mathbf{D}+%
\mathbf{D}^{\star }\mathbf{D})(t)\in D^{t}$ is an $\mathcal{A}_{t}^{\prime }$
linear combination of $\{\mathbf{D}_{i}(t)\}$ and a $\mathbf{G}\in D_{0}^{t}$%
, as well as $b^{\ast }b\in \mathcal{C}$ for $b=y_{0}+y_{i}\lambda ^{i}\in 
\mathcal{C}$ is a linear combination of $\widehat{y}^{i}$ and a $\widehat{b}%
\in \mathcal{C}_{0}$. Hence, $B^{\ast }B$ is a process of the same form as $%
B(t)$, what means that the operators $B(t)$ compose a $^{\ast }$-subalgebra $%
\mathcal{C}^{t}$ of $\mathcal{A}_{t}^{\prime }$. The algebra $b$ is a weakly
dense in $\mathcal{A}_{t}^{\prime }$ because it has the same commutant $%
\mathcal{A}_{t}$, as the family $\{Y(s):s\geq t\}$, and hence generates the
same von Neumann algebra $\mathcal{A}_{t}^{\prime }$. \medskip

\noindent \textsc{Proof of the Theorem \ref{th:4}.\/} We shall look for the
martingale $M$, defining the decomposition (\ref{eq:4.8}) in the Lemma \ref%
{lemma:4.1}. Let us suppose, that it is a stochastic integral
nonanticipating span 
\begin{equation*}
M_{t}=\int_{0}^{t}(\mathbf{Z}^{\star }\mathbf{\hat{D}}_{i}\mathbf{Z})_{\nu
}^{\mu }(s)\kappa _{s}^{i}\mathrm{d}A_{\mu }^{\nu }(s)\ ,\quad \kappa
_{t}^{i}\in \mathcal{A}_{t}^{1}
\end{equation*}%
of the observable martingales 
\begin{equation*}
\tilde{Y}_{i}(t)=Y_{i}(t)-\int_{0}^{t}\epsilon _{s}(\mathbf{Z}^{\star }%
\mathbf{D}_{i}\mathbf{Z})_{+}^{-}(s)\mathrm{d}s=\int \mathbf{Z}^{\star }%
\mathbf{\hat{D}}_{i}\mathbf{Z}\mathrm{d}\widehat{A}\ ,
\end{equation*}%
where $\tilde{Y}_{i}(t)$ should not be taken into account, if $D_{i}(t)\in 
\mathcal{D}_{0}^{t}$, as it is a zero almost everywhere martingale according
to the Lemma \ref{lemma:4.2}. Due to the weak density of $\mathcal{C}^{t}$
in $\mathcal{A}_{t}^{\prime }$, proved in Lemma \ref{lemma:4.3}, it is
sufficient to find the coefficient $\kappa _{t}^{i}$ from the condition 
\begin{equation*}
<\xi |B(t)^{\ast }X(t)\xi >=<\xi |B(t)^{\ast }\epsilon _{t}(X(t))\xi >
\end{equation*}%
for all $B(t)$ in the form (\ref{eq:4.11}). Using the QS It\^{o} formula (%
\ref{eq:2.4}) for $\mathrm{d}B=(\mathbf{Z}^{\star }\mathbf{DZ})_{\nu }^{\mu }%
\mathrm{d}A_{\mu }^{\nu }$, where $\mathbf{D}=\mathbf{D}_{0}+\mathbf{\hat{D}}%
_{i}\lambda ^{i}$, one can obtain 
\begin{eqnarray*}
\mathrm{d} &<&\xi |B^{\ast }X\xi >=<\xi |(\mathbf{Z}^{\star }(B\otimes 
\mathbf{\delta }+\mathbf{D})^{\star }\mathbf{FZ})_{+}^{-}\xi >\mathrm{d}t= \\
&=&<\xi |B^{\ast }\epsilon _{t}(\mathbf{Z}^{\star }\mathbf{FZ})_{+}^{-}+(%
\mathbf{Z}^{\star }\mathbf{D}^{\star }\mathbf{FZ})_{+}^{-}|\xi >\mathrm{d}t\
.
\end{eqnarray*}%
On the other hand, taking into account that $\mathrm{d}\epsilon
_{t}(X(t))=\epsilon _{t}(\mathbf{Z}^{\star }\mathbf{FZ})_{+}^{-}(t)\mathrm{d}%
t+\mathrm{d}M_{t}$, $\mathrm{d}M=(\mathbf{Z}^{\star }\mathbf{\hat{D}}_{i}%
\mathbf{Z})_{\nu }^{\mu }\kappa ^{i}\mathrm{d}\widehat{A}_{\mu }^{\nu }$,
one can obtain 
\begin{eqnarray*}
\mathrm{d} &<&\xi |B^{\ast }\epsilon _{t}(X)\xi >=<\xi |B^{\ast }\epsilon
_{t}(\mathbf{Z}^{\star }\mathbf{FZ})_{+}^{-}\xi >\mathrm{d}t+ \\
+ &<&\xi |(\mathbf{Z}^{\star }\mathbf{D}^{\star }\mathbf{Z})_{+}^{-}\epsilon
_{t}(X)+(\mathbf{Z}^{\star }\mathbf{D}^{\star }\widetilde{\mathbf{D}}_{i}%
\mathbf{Z})_{+}^{-}\kappa _{t}^{i}|\xi >\mathrm{d}t\ .
\end{eqnarray*}%
Hence, $<\xi |\mathbf{Z}^{\star }\mathbf{D}^{\star }\widetilde{\mathbf{D}}%
_{i}\mathbf{Z}|\kappa _{t}^{i}\xi >_{+}^{-}=<\xi |\left\{ \left( \mathbf{Z}%
^{\star }\mathbf{D}^{\star }\mathbf{FZ}\right) _{+}^{-}-\left( \mathbf{Z}%
^{\star }\mathbf{D}^{\star }\mathbf{Z}\right) _{+}^{-}\epsilon
_{t}(X)\right\} \xi >$, what is equivalent for $\mathbf{D}=\mathbf{D}_{0}+%
\widetilde{\mathbf{D}}_{i}\lambda ^{i}$ to (\ref{eq:4.5}) and 
\begin{equation}
<\xi |(\mathbf{Z}^{\star }\mathbf{D}_{0}^{\star }\mathbf{D}_{i}\mathbf{Z}%
)_{+}^{-}\kappa _{t}^{i}\xi >=<\xi |\{(\mathbf{Z}^{\star }\mathbf{D}%
_{0}^{\star }\mathbf{FZ})_{+}^{-}-(\mathbf{Z}^{\star }\mathbf{D}_{0}^{\star }%
\mathbf{Z})_{+}^{-}\epsilon _{t}(X)\}\xi >  \label{eq:4.12}
\end{equation}%
due to $\mathbf{D}^{\star }\widetilde{\mathbf{D}}_{i}=\mathbf{D}^{\star }%
\mathbf{D}_{i}$ and arbitrariness of $\lambda _{t}^{i}\in \mathcal{A}%
_{t}^{\prime }$. But the left hand side of the last equation (\ref{eq:4.12})
due to Schwarz inequality is zero: 
\begin{equation*}
<\xi |(\mathbf{Z}^{\star }\mathbf{D}_{0}^{\star }\mathbf{D}_{0}\mathbf{Z}%
)_{+}^{-}\xi >=<(D_{\circ }^{\circ }Z_{+}^{\circ }+D_{+}^{\circ })_{0}\xi
|(D_{\circ }^{\circ }Z_{+}^{\circ }+D_{+}^{\circ })_{i}\xi >=0
\end{equation*}%
as $<\xi |(\mathbf{Z}^{\star }\mathbf{D}_{0}^{\star }\mathbf{D}_{0}\mathbf{Z}%
)_{+}^{-}\xi >=\Vert (D_{\circ }^{\circ }Z_{+}^{\circ }+D_{+}^{\circ
})_{0}\xi \Vert ^{2}=0$ for $D_{0}\in \mathcal{D}_{0}^{t}$. On the other
hand, taking into account, that 
\begin{equation*}
<\xi |(\mathbf{Z}^{\star }\mathbf{D}_{0}^{\star }\mathbf{FZ})_{+}^{-}\xi
>=<\xi |(\mathbf{Z}^{\star }\mathbf{D}_{0}(X\otimes \mathbf{\delta })\mathbf{%
Z})_{+}^{-}\xi >,
\end{equation*}%
as $<\xi |(\mathbf{Z}^{\star }\mathbf{D}_{0}^{\star }\mathbf{CZ})_{+}^{-}\xi
>=<(D_{\circ }^{\circ }Z_{+}^{\circ }+D_{+}^{\circ })_{0}\xi |(C_{\circ
}^{\circ }Z_{+}^{\circ }+C_{+}^{\circ })\xi >=0$ for $\mathbf{C}=\mathbf{F}%
-X\otimes \mathbf{\delta }$, and due to $\star $ -- normality 
\begin{equation*}
D_{\circ }^{\circ }D_{\circ }^{\circ \ast }=D_{\circ }^{\circ \ast }D_{\circ
}^{\circ },\;D_{\circ }^{\circ }D_{\circ }^{-\ast }=D_{\circ }^{\circ \ast
}D_{+}^{\circ },\;D_{\circ }^{-}D_{\circ }^{-\ast }=D_{+}^{\circ \ast
}D_{+}^{\circ }
\end{equation*}%
of $\mathbf{D}_{0}\in \mathcal{D}_{0}^{t}$ as for a matrix--operator of the
commutative matrix $\star $--algebra ${\mathcal{F}}_{t}^{\prime }$, one can
obtain 
\begin{eqnarray*}
\lefteqn{<\xi |(Z^{\star }\mathbf{D}_{0}^{\star }FZ)_{+}^{-}\xi >=<\xi
|(D_{+}^{-}+D_{\circ }^{-}Z_{+}^{\circ })_{0}^{\ast }X\xi >+} \\
&+<&(D_{\circ }^{\circ }Z_{+}^{\circ }+D_{+}^{\circ })_{0}\xi |XZ_{+}^{\circ
}\xi >=<\xi |(D_{+}^{-\ast }X+Z_{+}^{\circ \ast }XD_{\circ }^{\circ \ast
}Z_{+}^{\circ })_{0}\xi >= \\
&=&<\xi |(D_{+}^{-\ast }+D_{+}^{\circ \ast }D_{\circ }^{\circ +}D_{+}^{\circ
})_{0}X\xi >=<\xi |(D_{+}^{-}+D_{+}^{\circ \ast }D_{\circ }^{\circ
+}D_{+}^{\circ })_{0}^{\ast }\epsilon _{t}(X)\xi >.
\end{eqnarray*}%
Here $D_{\circ }^{\circ ^{+}}$ is quasi-inverse conjugate matrix--operator
for normal $D_{\circ }^{\circ }=D_{\circ }^{\circ ^{\ast }}D_{\circ }^{\circ
^{+}}D_{\circ }^{\circ }$, and we used 
\begin{equation*}
(D_{\circ }^{\circ }Z_{+}^{\circ }+D_{+}^{\circ })_{0}\xi =0,\quad (D_{\circ
}^{\circ \ast }Z_{+}^{\circ }+D_{\circ }^{-\ast })_{0}\xi =0
\end{equation*}%
as for $\star $ -- normal $\mathbf{D}_{0}\in \mathcal{D}_{0}^{t}$. Hence,
the right side of equation (\ref{eq:4.12}) is also zero: 
\begin{equation*}
<\xi |(\mathbf{Z}^{\star }\mathbf{D}_{0}^{\star }\mathbf{Z})_{+}^{-}\xi
>=<\xi |(\mathbf{Z}^{\star }\mathbf{D}_{0}^{\star }\mathbf{Z}%
)_{+}^{-}\epsilon _{t}(X)\xi >,
\end{equation*}%
because in the same way one can obtain 
\begin{eqnarray*}
<\xi |(\mathbf{Z}^{\star }\mathbf{D}_{0}^{\star }\mathbf{Z})_{+}^{-}\epsilon
_{t}(X)\xi > &=&<\xi |(D_{+}^{-}+D_{\circ }^{-}Z_{+}^{\circ })_{0}^{\ast
}\epsilon _{t}(X)\xi >+ \\
+<(D_{\circ }^{\circ }Z_{+}^{\circ }+D_{+}^{\circ })_{0}\xi |Z_{+}^{\circ
}\epsilon _{t}(X)\xi > &=&\xi |(D_{+}^{-}+D_{+}^{\circ }D_{\circ }^{\circ
+}D_{+}^{\circ })_{0}^{\ast }\epsilon _{t}(X)\xi >.
\end{eqnarray*}%
This proves also the uniqueness of the solution of the equation (\ref{eq:4.5}%
) up to the kernel of the correlation matrix $\varrho (t)=\left( \varrho
_{ik}\right) (t)$ because if $\varrho _{ik}\lambda _{0}^{k}=0$ for an $%
\mathcal{A}_{t}^{\prime }$ -- adapted vector process $\lambda _{0}=(\lambda
_{0}^{i})$, then $\mathbf{D}_{0}=\lambda _{0}^{i}\mathbf{D}_{i}\in \mathcal{D%
}_{0}$, and, hence, 
\begin{eqnarray*}
<\lambda ^{i}\xi |\left\{ \epsilon (\mathbf{Z}^{\star }\mathbf{D}_{i}^{\star
}\mathbf{FZ})_{+}^{-}-\epsilon (\mathbf{Z}^{\star }\mathbf{D}_{i}^{\star }%
\mathbf{Z})_{+}^{-}\epsilon (X)\right\} \xi > &=& \\
+<\xi |\left\{ (\mathbf{Z}^{\star }\mathbf{D}_{0}^{\star }\mathbf{FZ}%
)_{+}^{-}(\mathbf{Z}^{\star }\mathbf{D}_{0}^{\star }\mathbf{Z}%
)_{+}^{-}\epsilon (X)\right\} \xi > &=&0.
\end{eqnarray*}%
In the case $F_{\nu }^{\mu }=X\delta _{\nu }^{\mu }$ taking into account
that 
\begin{equation*}
\epsilon (\mathbf{Z}^{\star }\mathbf{D}(X\otimes \mathbf{\delta })\mathbf{Z}%
)_{+}^{-}=D_{+}^{-}\epsilon (X)+D_{\circ }^{-}\epsilon (XZ_{+}^{\circ
})+\epsilon (Z_{+}^{\circ \ast }X)D_{+}^{\circ }+\epsilon (Z_{+}^{\circ \ast
}D_{\circ }^{\circ }XZ_{+}^{\circ })
\end{equation*}%
and $\epsilon (\mathbf{Z}^{\star }\mathbf{DZ})_{+}^{-}\epsilon (Z)=$ 
\begin{equation*}
D_{+}^{-}\epsilon (X)+D_{\circ }^{-}\epsilon (X)\epsilon (Z_{+}^{\circ
})+\epsilon (Z_{+}^{\circ })^{\ast }\epsilon (X)D_{+}^{\circ }+\epsilon
(Z_{+}^{\circ \ast }D_{\circ }^{\circ }Z_{+}^{\circ })\epsilon (X),
\end{equation*}%
one can easily obtain the equation (\ref{eq:4.6}) from (\ref{eq:4.5}).

One should look for the initial condition $\epsilon _{0}(\widehat{x}\otimes
I)=\epsilon (\widehat{x})\otimes I$ for the equation (\ref{eq:4.2}) in the
linear form $\epsilon (\widehat{x})=<\psi |\widehat{x}\psi >+\tilde{y}%
_{i}\kappa ^{i}$, where $\kappa ^{i}$ should be found from $<\widehat{b}\psi
|\widehat{z}\psi >=<\widehat{b}\psi |\epsilon (\widehat{z})\psi >$ for all $%
\widehat{b}=\widehat{y}_{0}+\Sigma \tilde{y}_{i}\lambda ^{i}$, where $%
\widehat{y}_{0}\in \mathcal{C}$ and $\lambda ^{i}\in \mathbb{C}$. This gives
the initial equation (\ref{eq:4.7}).\hfill \vrule height .9ex width .8ex
depth -.1ex 

\begin{corollary}
If $\{y_{i}\}$ are commuting orthogonal projectors in $\mathfrak{h}$, and
also $\mathbf{D}_{i}^{\star }=\mathbf{D}_{i}=\mathbf{D}_{i}^{2}$ are
commuting $\star $-matrix projectors, then the conditions of the Theorem \ref%
{th:4} are fullfilled, and they are fullfilled also in the case $D_{\circ
i}^{\circ }=0$. In particular, for the case $D_{\circ }^{\circ }=I\otimes
\delta _{\circ }^{\circ }$, $D_{\circ }^{-}=0=D_{+}^{\circ }$, $D_{+}^{-}=0$%
, corresponding to the counting output process $Y=N$, the equation (\ref%
{eq:4.6}) gives 
\begin{equation*}
\kappa _{t}=\epsilon _{t}(Z_{+}^{\circ \ast }XZ_{+}^{\circ })(t)/\epsilon
_{t}(Z_{+}^{\circ \ast }Z_{+}^{\circ })(t)-\epsilon _{t}(X)(t)\ ,
\end{equation*}%
%
%
%
%
if $\epsilon _{t}(Z_{+}^{\circ \ast }Z_{+}^{\circ })(t)\not=0$. In the other
case $D_{\circ }^{-}=1=D_{+}^{\circ }$, $D_{\circ }^{\circ }=0$, $%
D_{+}^{-}=0 $, corresponding to the output coordinate observation $Y=Q$, one
obtains 
\begin{equation}
\kappa _{t}=\epsilon _{t}(XZ_{+}^{\circ }+Z_{+}^{\circ \ast }X)-\epsilon
_{t}(X)\epsilon _{t}(Z_{+}^{\circ }+Z_{+}^{\circ \ast })\ .  \label{eq:4.13}
\end{equation}%
\label{corol:3}
\end{corollary}

Indeed, the linear, span of commuting orthoprojectors $\{y_{i}\}$, and also $%
\mathcal{A}_{t}^{\prime }$-span of $\star $-projectors $\{D_{i}\}$ is a $%
\ast $-and $\star $-algebra $\mathcal{C}$ and $\mathcal{D}^{t}$
correspondingly. In the case $D_{i\circ }^{\;\circ }=0$ the product $\mathbf{%
D}_{i}^{\star }\cdot \mathbf{D}_{k}$ is in $\mathcal{D}_{0}^{t}$ as matrix
with $(\mathbf{D}_{i}^{\star }\mathbf{D}_{k})_{\circ }^{\circ }=0$, $(%
\mathbf{D}_{i}^{\star }\mathbf{D}_{k})_{+}^{\circ }=0$, $(\mathbf{D}%
_{i}^{\star }\mathbf{D}_{k})_{\circ }^{-}=0$. Such commutative matrices also
form a commutative $\star $-algebra $\mathcal{D}^{t}$ up to the ideal $%
\mathcal{D}_{0}^{t}$, because $\mathbf{G}^{\star }\mathbf{G}=0$, and, hence, 
$(\mathbf{GZ}_{+}|\mathbf{GZ}_{+})=0$ for a matrix-operator $\mathbf{G}$
with $G_{\nu }^{\mu }=0$ for $(\mu ,\nu )\not=(-,+)$. 

\section{An application of the QS filtering}

\label{sec:belat5}

The applications of the filtering equation (\ref{eq:4.2}) to the derivation
of a posteriori Schr\"odinger equation for the coordinate observation are
given in \cite{bib:ref17,bib:ref18}, and for the counting observation are
given in \cite{bib:ref19,bib:ref21}.

In contrast to the usual Schr\"{o}dinger equation, describing a closed
quantum system without observation, these new stochastic wave equations give
the dynamics of an open quantum system undergoing the nondemolition
measurements which are continuous in time. Thus the continual wave packed
reduction problem is solved by the quantum filtering method for the typical
QS models of observation such as a quantum particle is a bubble chamber \cite%
{bib:ref18} (diffusive observation) and an atom radiating the photons \cite%
{bib:ref21} (counting observation). Here we consider another example of the
quantum nonlinear filtering -- the QS spin localization, describing the
continuous collapse of the vector polarization $\vec{p}=(p_{1},p_{2},p_{3})$
for the spin ${\frac{1}{2}}$ of an electron under a continuous nondemolition
measurement in a magnetic field. The polarization $\vec{p}(t)$ at the time
instant $t>0$ is given by the conditional expectations (\ref{eq:3.5}) 
\begin{equation}
p_{j}(t)=\epsilon _{t}(X_{j}(t))\ ,\ X_{j}(t)=U^{\ast }(t)(\hat{x}%
_{j}\otimes I)U(t)\ ,  \label{eq:5.1}
\end{equation}%
where $\hat{x}_{j}=\hat{\sigma}_{j}$ are the Pauli matrices 
\begin{equation*}
\hat{\sigma}_{1}=\left( 
\begin{array}{cc}
0 & 1 \\ 
1 & 0%
\end{array}%
\right) \ ,\quad \hat{\sigma}_{2}=\left( 
\begin{array}{cc}
0 & -i \\ 
i & 0%
\end{array}%
\right) \ ,\quad \hat{\sigma}_{3}=\left( 
\begin{array}{cc}
1 & 0 \\ 
0 & -1%
\end{array}%
\right)
\end{equation*}%
and $U(t)$ is a QS unitary evolution in the Hilbert space $\mathcal{H}%
^{0}\otimes \mathcal{F}$. Here $\mathcal{H}^{0}$ is the Hilbert space $%
\mathbb{C}^{2}\otimes L^{2}(\mathbb{R}^{3})$ of the spinors $\psi (\vec{r}%
)=\left( 
\begin{array}{c}
\psi _{-} \\ 
\psi _{+}%
\end{array}%
\right) (\vec{r})$, where $\psi _{\pm }(\vec{r})$, $\vec{r}\in \mathbb{R}%
^{3} $ are the wave functions of the nonrelativistic electron with the
definite $z $-projections $\pm {\frac{1}{2}}$ of its spin $\hat{\vec{s}}={%
\frac{1}{2}}(\hat{\sigma}_{1},\hat{\sigma}_{2},\hat{\sigma}_{3})$ and
probabilistic normalization $\Vert \psi \Vert ^{2}\int \psi (\vec{r}%
)^{+}\psi (\vec{r})\mathrm{d}r=1$, where $\psi ^{+}\psi =|_{-}|^{2}+|\psi
_{+}|^{2}$, and $\mathcal{F}=\Gamma (\mathcal{E})$ is the Fock space over
the Hilbert space $\mathcal{E}=L^{2}(\mathbb{R}^{+})\otimes \mathbb{C}^{n}$%
.\hfill \break The initial polarization $\vec{p}%
(0)=(p_{1}^{0},p_{2}^{0},p_{3}^{0})=\vec{p}_{0}$, 
\begin{equation*}
p_{j}^{0}=\int \psi (\vec{r})^{\dag }\hat{\sigma}_{j}\psi (\vec{r})\mathrm{d}%
\vec{r}\ ,\quad j=1,2,3\ ,
\end{equation*}%
has the values in the unite ball $\mathcal{B}=\{\vec{p}\in \mathbb{R}^{3}:\;|%
\vec{p}|\leq 1\}$, where $|\vec{p}|^{2}=(\vec{p},\vec{p})\equiv p^{2}$, i.e.
can be mixed $\sum_{j=1}^{3}(p_{j}^{0})^{2}<1$ even in the pure (vector)
state $\psi \in \mathcal{H}^{0}$, $\Vert \psi \Vert ^{2}=\int \psi (\vec{r}%
)^{\dag }\psi (\vec{r})\mathrm{d}\vec{r}=1$.

Let us suppose that the evolution $U(t)$ defines the system of Langevin
equations (\ref{eq:2.14}) of the form 
\begin{equation}
\mathrm{d}\vec{X}+(\mathrm{i}[\vec{X},H]+{\frac{1}{2}}\sum_{i=1}^{n}[[\vec{X}%
,L_{j}],L_{j}])\mathrm{d}t=\mathrm{i}\sum_{j=1}^{n}[\vec{X},L_{i}]\mathrm{d}%
V_{i}\ .  \label{eq:5.2}
\end{equation}%
Here $\vec{X}(t)=(X_{1},X_{2},X_{3})(t)$, $H(t)={\frac{1}{2}}%
\sum_{j=1}^{3}u^{j}(t)X_{j}(t)$ is the spin-Hamiltonian, corresponding to
the magnetic tense $\vec{u}(t)=(u^{1},u^{2},u^{3})(t)\in \mathbb{R}^{3}$,%
\begin{equation*}
L_{i}(t)={\frac{1}{2}}\sum_{j=1}^{n}r_{i}^{j}(t)X_{j}(t)\equiv {\frac{1}{2}}%
R_{i}(t)
\end{equation*}%
are spin-operators, defined by the real vectors $\vec{r}%
_{i}(t)=(r_{i}^{1},r_{i}^{2},r_{i}^{3})(t)\in \mathbb{R}^{3}$, $i=1,\dots ,n$%
, and $V_{i}=\hat{1}\otimes 2\Im A_{i}^{+}$ $i=1,\dots ,n$ are the
independent standard Wiener processes, represented by the input operators ${%
\frac{1}{\mathrm{i}}}(A_{i}^{+}(t)-A_{-}^{i}(t))$, $\mathrm{i}=\sqrt{-1}$ in
the Fock space $\mathcal{F}$ with respect to the initial vacuum state $%
\delta _{\emptyset }\in \mathcal{F}$. The stochastic system of the operator
equations (\ref{eq:5.2}) corresponds to the unitary Markovian evolution (\ref%
{eq:2.6}) in $\mathfrak{h}\otimes \mathcal{F}$, $\mathfrak{h}=\mathbb{C}^{2}$
with the generators 
\begin{equation*}
Z_{\nu }^{\mu }(t)=U(t)^{\ast }(\hat{Z}_{\nu }^{\mu }(t)\otimes I)U(t)%
\begin{array}{c}
\mu =-,i,\dots ,n \\ 
\nu =1,\dots ,n,+%
\end{array}%
\end{equation*}%
defined by the spin-operators 
\begin{equation*}
\hat{z}_{k}^{i}=\delta _{k}^{i}\hat{1}\ ,\quad \hat{z}_{+}^{i}={\frac{1}{2}}%
\hat{r}_{i}=\hat{z}_{i}^{-}\ ,\ \hat{z}_{+}^{-}=-{\frac{1}{2}}\left( {\frac{1%
}{4}}\hat{r}^{2}+\mathrm{i}\hat{u}\right) \ ,\quad \mathrm{i}=\sqrt{-1}\ ,
\end{equation*}%
where $\hat{r}_{i}(t)=\sum_{i=1}^{3}r_{i}^{j}(t)\hat{\sigma}_{j}$, $\hat{r}%
^{2}(t)=\sum_{j=1}^{n}r_{i}^{2}(t)\hat{1}$, $\hat{u}(t)=\sum%
\limits_{j=1}^{3}u^{j}(t)\hat{\sigma}_{j}$. Such the evolution realizes the
output coordinate processes 
\begin{equation*}
Y_{i}(t)=U(t)^{\ast }W_{i}(t)U(t)=Q_{i}(t),\quad i=1,\dots ,n,
\end{equation*}%
satisfying the QS equations (\ref{eq:2.5}) in the form 
\begin{equation}
\mathrm{d}Y_{i}=R_{i}\mathrm{d}t+\mathrm{d}W_{i}\ ,\quad \mathrm{d}W_{i}=%
\hat{1}\otimes 2\Re A_{i}^{+}  \label{eq:5.3}
\end{equation}%
of the indirect non-demolition observation of the noncommuting
spin-operators 
\begin{equation*}
R_{i}(t)=U(t)^{\ast }(\hat{r}_{i}(t)\otimes I)U(t)\ ,\quad i=1,\dots ,n\ .
\end{equation*}

The standard Wiener processes $W_{i}$, $i=1,\dots ,n$, represented by the
commuting operators $A_{i}^{+}(t)+A_{-}^{i}(t)$ in $\mathcal{F}$, describe
the independent errors $\dot{Y}_{i}-R_{i}$ as the white noises $\dot{W}_{i}$%
. They do not commute with the white noises $\dot{V}_{i}$ of the
perturbations in the quantum system (\ref{eq:5.2}): 
\begin{equation}
\lbrack \dot{V}_{i}(s)\ ,\ \dot{W}_{k}(t)]=2\mathrm{i}\delta (s-t)\delta
_{ik}\hat{I}\ ,  \label{eq:5.4}
\end{equation}%
due to 
\begin{equation*}
\lbrack V_{i}(s)\ ,\ \dot{W}_{k}(t)]=2\mathrm{i}[\hat{A}_{-}^{i}(s),\ \hat{A}%
_{k}^{+}(t)]=2\mathrm{i}\min (s,t)\delta _{ik}\hat{I}\ .
\end{equation*}%
%
%
%
%

\begin{proposition}
Under the given assumptions the a posteriori spin polarizations (\ref{eq:5.1}%
) satisfy the following system of nonlinear stochastic equations 
\begin{equation}
\mathrm{d}\vec{p}+(\vec{p}\wedge \vec{u}+{\frac{1}{2}}%
\sum_{i=1}^{n}(r_{i}^{2}\vec{p}-(\vec{p},\vec{r}_{i})\vec{r}_{i}))\mathrm{d}%
t=\sum_{i=1}^{n}(\vec{r}_{i}-(\vec{p},\vec{r}_{i})\vec{p})\mathrm{d}\tilde{Y}%
_{i}\ ,  \label{eq:5.5}
\end{equation}%
where $\mathrm{d}\tilde{Y}_{i}(t)=\mathrm{d}Y_{i}(t)-(\vec{p}(t)$, $\vec{r}%
_{i}(t))\mathrm{d}t$. \label{prop:5}
\end{proposition}

\noindent \textsc{Proof.\/} Let us consider the nonlinear filtering equation
(\ref{eq:4.2}) for the spin-operators $X_{j}(t)$, $j=1,2,3$, which are
equivalent to the Pauli matrices $\hat{\sigma}_{j}$, $j=1,2,3$. We can use (%
\ref{eq:4.2}) for the evaluation of the expectations (\ref{eq:5.1}) because
the conditions of the Theorem \ref{th:4} are fulfilled (see the Corollary of
sec. \ref{sec:belat4}). The innovating martingales $\mathrm{d}\tilde{Y}_{i}=%
\mathrm{d}Y_{i}-\epsilon _{t}(r_{i})\mathrm{d}t$ in this case are given by
the differences $\mathrm{d}\tilde{Y}_{i}=\mathrm{d}Y_{i}-(\vec{p},\vec{r}%
_{i})\mathrm{d}t$ because 
\begin{equation*}
\epsilon _{t}(R_{i}(t))=\sum_{j=1}^{3}r_{i}^{j}(t)\epsilon
_{t}(X_{j}(t))=\sum_{j=1}^{3}r_{i}^{j}(t)p_{j}(t)\ .
\end{equation*}%
Due to $\rho _{ik}(t)\mathrm{d}t=\epsilon _{t}(\mathrm{d}\tilde{Y}_{i}(t)%
\mathrm{d}\tilde{Y}_{k}(t))=\delta _{ik}\mathrm{d}t$, the coefficients $%
\kappa _{t}^{i}(X_{j}(t))$ are given by 
\begin{eqnarray*}
\kappa _{t}^{i}(\vec{x}) &=&{\frac{1}{2}}\epsilon _{t}(\vec{X}%
(t)R_{i}(t)+R_{i}(t)\vec{X}(t))-\epsilon _{t}(\vec{X}(t))\epsilon
_{t}(R_{i}(t))= \\
&=&\vec{r}_{i}(t)-(\vec{p}(t),\ \vec{r}_{i}(t))\vec{p}(t)\ ,
\end{eqnarray*}%
because $\hat{\sigma}_{j}\hat{r}_{j}+\hat{r}_{j}\hat{\sigma}_{j}=2r_{i}^{j}%
\hat{1}$ for $\hat{r}_{i}=\sum\limits_{j=1}^{3}r_{i}^{j}\hat{\sigma}_{j}$,
and 
\begin{equation*}
X_{j}R_{i}+R_{i}X_{j}=U(t)^{\ast }(\hat{\sigma}_{j}\hat{r}_{i}+\hat{r}_{i}%
\hat{\sigma}_{j})U(t)=2r_{i}^{j}\hat{I}\ .
\end{equation*}

The vector-product $\vec{p}(t)\wedge \vec{u}(t)$ in (\ref{eq:5.5})
represents the expectations 
\begin{equation*}
\mathrm{i}\epsilon _{t}([\vec{X}(t)\ ,\ H(t)])=\mathrm{i}\epsilon _{t}\left( 
{\frac{1}{2}}\sum_{j=1}^{3}[\vec{X}(t)\ ,\ X_{j}(t)]u_{j}(t)\right) \ ,
\end{equation*}%
because $[\hat{\vec{\sigma}},\hat{u}]=\Sigma _{i=1}^{3}[\hat{\vec{\sigma}},%
\hat{\sigma}_{i}]u^{i}={\frac{2}{\mathrm{i}}}\hat{\vec{\sigma}}\wedge \vec{u}
$, and%
\begin{equation*}
\lbrack \vec{X}(t)\ ,\ H(t)]={\frac{1}{2}}U(t)^{\ast }([\hat{\vec{\sigma}},%
\hat{u}(t)]\otimes I)U(t).
\end{equation*}
In the same way one can obtain 
\begin{equation*}
r_{i}^{2}(t)\vec{p}(t)-(\vec{p}(t),\vec{r}_{i}(t))\vec{r}_{i}(t)=(\vec{p}%
(t)\wedge \vec{r}_{i}(t))\wedge \vec{r}_{i}(t)
\end{equation*}%
as for the vector representation of the double commutator ${\frac{1}{2}}[[%
\vec{X}(t)\ ,\ L_{i}(t)]\ ,\ L_{i}(t)]$, defining together with $\mathrm{i}[%
\vec{X}(t)\ ,\ H(t)]$ the products $(\mathbf{Z}^{\star }(t)\vec{X}(t)\mathbf{%
Z}(t))_{+}^{-}$ in (\ref{eq:4.2}).\hfill \vrule height .9ex width .8ex depth
-.1ex

Now we can prove, that the continuous indirect nondemolition measurement (%
\ref{eq:5.3}) of the quantum spin reduces any initial state of the electron
at the limit $t\to\infty$, to the completely polarized one. This gives a
kind of the stochastic ergodicity property of the nonlinear system of
quantum filtering equation (\ref{eq:5.5}). 

\begin{theorem}
Let $\vec{p}(0)=\vec{p}_{0}\in \mathcal{B}$ be an arbitrary initial
polarization for the nonlinear quantum filtering equation (\ref{eq:5.5}).
Then this equation has a unique stochastic solution $\vec{p}(t)\in \mathcal{B%
}$, and $p^{2}(t)=(\vec{p}(t)\ ,\ \vec{p}(t))\rightarrow 1$ at $t\rightarrow
\infty $ almost surely, if $\lambda
(t)=\int_{0}^{t}\sum\limits_{i=1}^{n}|r_{i}(s)|^{2}\mathrm{d}s\rightarrow
\infty $. \label{th:5}
\end{theorem}

\noindent \textsc{Proof.\/} The vector stochastic equation (\ref{eq:5.5}) up
to a renormalization $\vec{f}(t)=\rho (t)\vec{p}(t)$ is equivalent to the
linear stochastic equation 
\begin{equation}
\mathrm{d}\vec{f}+(\vec{f}\wedge \vec{u}+{\frac{1}{2}}%
\sum_{i=1}^{u}(r_{i}^{2}\vec{f}-(\vec{f},\vec{r}_{i})\vec{r}_{i}))\mathrm{d}%
t=\rho \sum_{i=1}^{n}\vec{r}_{i}\mathrm{d}Y_{i}\ .  \label{eq:5.6}
\end{equation}%
Indeed, let $\rho (t)$ be the stochastic It\^{o}'s integral 
\begin{equation}
\rho (t)=1+\int_{0}^{t}\sum_{i=1}^{u}(\vec{f}(s),\ \vec{r}_{i}(s))\mathrm{d}%
Y_{i}(s)  \label{eq:5.7}
\end{equation}%
defined by the unique solution $\vec{f}(t)$ of this ordinary linear
stochastic differential equation with the initial nonstochastic vector $\vec{%
f}(0)=\vec{p}_{0}$. Then $\mathrm{d}\rho =\sum\limits_{i=1}^{u}(\vec{f},\vec{%
r}_{i})\mathrm{d}Y_{i}$, and by It\^{o}'s formula 
\begin{equation*}
\mathrm{d}(\rho \vec{p})=\mathrm{d}\rho \vec{p}+\mathrm{d}\rho \mathrm{d}%
\vec{p}+\rho \mathrm{d}\vec{p}
\end{equation*}%
we obtain the equation for $\vec{f}=\rho \vec{p}$ iff $\vec{p}(t)$ satisfies
the equation (\ref{eq:5.5}): 
\begin{eqnarray*}
\lefteqn{\mathrm{d}\vec{f}+(\vec{f}\wedge \vec{u}+\frac{1}{2}%
\sum_{i=1}^{u}r_{i}^{2}\vec{f}-(\vec{f},\vec{r}_{i})\vec{r}_{i}))\mathrm{d}t=%
\mathrm{d}\rho \vec{p}+\mathrm{d}\rho \mathrm{d}\vec{p}+} \\
&+&\rho \sum_{i=1}^{u}(\vec{r}_{i}-(\vec{p},\vec{r}_{i})\vec{p})\mathrm{d}%
\tilde{Y}_{i}=\sum_{i=1}^{u}(\vec{f},\vec{r}_{i})\vec{p}\mathrm{d}%
Y_{i}+\sum_{i=1}^{n}(\vec{r}_{i}-(\vec{p},\vec{r}_{i})\vec{p})(\vec{f},\vec{r%
}_{i})\mathrm{d}t \\
&+&\sum_{i=1}^{n}(\rho \vec{r}_{i}-(\vec{f},\vec{r}_{i})\vec{p})(\mathrm{d}%
Y_{i}-(\vec{p},\vec{r}_{i})\mathrm{d}t)=\rho \sum_{i=1}^{n}\vec{r}_{i}%
\mathrm{d}Y_{i}\ .
\end{eqnarray*}%
This the unique solution of the nonlinear filtering equation (\ref{eq:5.5})
with $\vec{p}(0)=\vec{p}_{0}$ can be written almost surely $(\rho (t)\not=0)$
as $\vec{p}(t)=\vec{f}(t)/\rho (t)$, where $\vec{f}(t)$ is the solution of
the linear equation (\ref{eq:5.6}) with $\vec{f}(0)=\vec{p}_{0}$, and $\rho
(t)$ is the integral (\ref{eq:5.7}).

In order to prove that almost surely $|\vec{p}(t)|\leq 1$, if $|\vec{p}%
_{0}|\leq 1$, it is sufficient to show, that 
\begin{equation*}
f^{2}(t)=(\vec{f}(t),\vec{f}(t))\leq \rho (t)^{2}\ \quad \mathrm{if}\quad 
\vec{f}(0)=\vec{p}_{0}\ .
\end{equation*}%
Using the It\^{o}'s formula we obtain 
\begin{eqnarray*}
\mathrm{d}f^{2} &=&2(\vec{f},\ \mathrm{d}\vec{f})+(\mathrm{d}\vec{f}\ ,\ 
\mathrm{d}\vec{f})=2\rho \sum_{i=1}^{u}(\vec{f},\vec{r}_{i})\mathrm{d}Y_{i}-
\\
&-&\sum_{i=1}^{u}(r_{i}^{2}f^{2}-(\vec{f},\vec{r}_{i})^{2}-\rho
^{2}r_{i}^{2})\mathrm{d}t=\mathrm{d}\rho ^{2}+(\rho
^{2}-f^{2})\sum_{i=1}^{u}r_{i}^{2}\mathrm{d}t\ ,
\end{eqnarray*}%
where $\mathrm{d}\rho ^{2}=2\rho \mathrm{d}\rho +(\mathrm{d}\rho )^{2}=2\rho
\sum\limits_{i=1}^{u}(\vec{f}\ ,\ \vec{r}_{i})(\mathrm{d}Y_{i}+(\vec{f},\vec{%
r}_{i}))$. Hence%
\begin{equation*}
\mathrm{d}(f^{2}-\rho ^{2})=\dot{\lambda}(\rho ^{2}-f^{2})\mathrm{d}t,
\end{equation*}
where $\dot{\lambda}=\sum\limits_{i=1}^{u}r_{i}^{2}\geq 0$, and 
\begin{equation*}
\rho ^{2}(t)-(\vec{f}(t),\ \vec{f}(t))=e^{-\lambda (t)}(1-(\vec{p}_{0},\ 
\vec{p}_{0})),\quad \forall t\ .
\end{equation*}%
Thus $f^{2}(t)\leq \rho ^{2}(t)$, if $|\rho ^{0}|\leq 1$, and $%
f^{2}(t)\rightarrow \rho ^{2}(t)$ exponentially at $t\rightarrow \infty $,
if $\lambda (t)\rightarrow \infty $ ($f^{2}(t)=\rho ^{2}(t)\ ,\ \forall t$,
if $|\vec{p}_{0}|=1$). This proves that $\vec{p}(t)=\vec{f}(t)/\rho
(t)\rightarrow 1$ almost surely $(\rho (t)\not=0)$ due to the positivity of $%
\rho (t)$.\hfill \vrule height .9ex width .8ex depth -.1ex 

\noindent \textbf{Remark.\/} The model (\ref{eq:5.2}) of continual
nondemolition measurements of noncommuting spin-operators $R_{i}(t)$, $%
i=1,\dots ,n$ in the quantum stochastic system (\ref{eq:5.2}) is unique in
the Fock space $\mathcal{F}=\Gamma (\mathcal{E})$ over the minimal Hilbert
space $\mathcal{E}=L^{2}(\mathbb{R}_{+})\otimes \mathbb{C}^{n}$. It can not
be realized in the framework of classical probability theory due to the
noncommutativity (\ref{eq:5.4}) of the quantum stochastic processes $%
V_{i}(t) $ and $W_{i}(t)$ though each of them can be described as the
classical one separately due to the selfnondemolition (commutativity)
property $[V_{i}(t)\ ,\ V_{k}(s)]=0=[W_{i}(t)\ ,\ W_{k}(s)]$.

The result obtained here in a rigorous mathematical way corresponds to a
rather intuitive physical picture of the continual spontaneous collapse of
the quantum spin under the non-demolition observation. This proves the
appropriateness of the given quantum stochastic setup for the theory of
continuous measurements and quantum filtering. 

\end{document}